# SECOND-ORDER ADJOINT SENSITIVITY ANALYSIS METHODOLOGY ($2^{nd}$-ASAM) FOR LARGE-SCALE NONLINEAR SYSTEMS:

# II. APPLICATION TO A NONLINEAR HEAT CONDUCTION BENCHMARK

Dan G. Cacuci

Department of Mechanical Engineering, University of South Carolina

E-mail: cacuci@cec.sc.edu

Corresponding author:

Department of Mechanical Engineering, University of South Carolina

300 Main Street, Columbia, SC 29208, USA

Email: cacuci@cec.sc.edu; Phone: (919) 909 9624;

**ABSTRACT**

This work presents an illustrative application of the *second-order adjoint sensitivity analysis methodology* ($2^{nd}$-*ASAM*) to a paradigm nonlinear heat conduction benchmark, which models a conceptual experimental test section containing heated rods immersed in liquid lead-bismuth eutectic. This benchmark admits an exact solution, thereby making transparent the underlying mathematical derivations. The general theory underlying $2^{nd}$-*ASAM* indicates that, for a physical system comprising $N_\alpha$ parameters, the computation of all of the first- and second-order response sensitivities requires (per response) at most $N_\alpha$ "large-scale" computations using the *first-level* and, respectively, *second-level adjoint sensitivity systems* ($1^{st}$-*LASS* and $2^{nd}$-*LASS*). For this illustrative problem, six "large-scale" adjoint computations sufficed to compute exactly all five $1^{st}$-order and fifteen distinct $2^{nd}$-order derivatives of the temperature response to the five model parameters. The construction and solution of the $2^{nd}$-*LASS* requires very little additional effort beyond the construction of the adjoint sensitivity system needed for computing the first-order sensitivities. Very significantly, only the *sources* on the right-sides of the heat conduction differential operator needed to be modified; the left-side of the differential equations (and hence the "solver" in large-scale practical applications) remains unchanged.

For the nonlinear heat conduction benchmark, the second-order sensitivities play the following roles:
   (a) They cause the "expected value of the response" to differ from the "computed nominal value of the response;" for the nonlinear heat conduction benchmark, however, these differences were insignificant over the range of temperatures (400-900K) considered.
   (b) They contribute to increasing the response variances and modifying the response covariances, but for the nonlinear heat conduction benchmark their contribution was smaller than that stemming from the $1^{st}$-order response sensitivities, over the range of temperatures (400-900K) considered.
   (c) They comprise the leading contributions to causing asymmetries in the response distribution. For the benchmark test section considered in this work, the heat source, the boundary heat flux, and the temperature at the bottom boundary of the test section would cause the temperature distribution in the test section to be skewed significantly towards values *lower* the the mean temperature. On the other hand, the model parameters entering the nonlinear, temperature-dependent, expression of the LBE conductivity would cause the temperature distribution in the test section to be skewed significantly towards values *higher* the the mean temperature. These opposite effects partially cancel each other. Consequently, the cumulative effects of model parameter uncertainties on the skewness of the temperature distribution in the test section is such that the temperature distribution in the LBE is skewed slightly toward higher temperatures in the cooler part of the test section, but becomes increasingly skewed towards temperatures lower than the mean temperature in the hotter part of the test section. Notably, the influence of the model parameter that controls the strength of the nonlinearity in the heat conduction coefficient for this LBE test section benchmark would be strong if it were the only uncertain model parameter. However, if all of the other model parameters are also uncertain, all having equal relative standard deviations, the uncertainties in the heat source and boundary heat flux diminish the impact of uncertainties in the nonlinear heat conduction coefficient for the range of temperatures (400-900K) considered for this LBE test section benchmark.

Ongoing work aims at generalizing the $2^{nd}$-ASAM to enable the exact and efficient computation of higher-order response sensitivities. The availability of such higher-order sensitivities is expected to affect significantly the fields of optimization and predictive modeling, including uncertainty quantification, data assimilation, model calibration and extrapolation.

# 1. INTRODUCTION

The accompanying PART I [1] of this work has presented the mathematical formalism of the "*<u>S</u>econd-<u>O</u>rder <u>A</u>djoint <u>S</u>ensitivity <u>A</u>nalysis <u>M</u>ethodology ($2^{nd}$-ASAM) for **nonlinear** systems.*" This is a new method for computing exactly and efficiently second-order sensitivities (i.e., functional derivatives) of nonlinear system responses (i.e., "system performance parameters" in physical, engineering, biological systems) to the parameters characterizing large-scale nonlinear systems. The definition of "system parameters" includes all computational input data, correlations, initial and/or boundary conditions, etc. The $2^{nd}$-ASAM builds on the *first-order adjoint sensitivity analysis methodology* ($1^{st}$-ASAM) for nonlinear systems originally introduced in [2, 3], and extends the work presented in [4]. For each functional-type response of interest in a physical system comprising $N_\alpha$ parameters and $N_r$ responses, the $2^{nd}$-ASAM requires one large-scale computation using the *first-level adjoint sensitivity system* ($1^{st}$-LASS) for obtaining all of the first-order sensitivities $\{\partial R(\mathbf{u},\boldsymbol{\alpha})/\partial\alpha_i\}_{(\mathbf{u}^0,\boldsymbol{\alpha}^0)}$, $i=1,...,N_\alpha$, followed by at most $N_\alpha$ large-scale computations using the *second-level adjoint sensitivity systems* ($2^{nd}$-LASS) for obtaining exactly all of the second-order sensitivities $\{\partial^2 R(\mathbf{u},\boldsymbol{\alpha})/\partial\alpha_j\partial\alpha_i\}_{(\mathbf{u}^0,\boldsymbol{\alpha}^0)}$, $i,j=1,...,N_\alpha$. In practice, however, the number of large-scale computations required for computing exactly all of the second-order sensitivities $\{\partial^2 R(\mathbf{u},\boldsymbol{\alpha})/\partial\alpha_j\partial\alpha_i\}_{(\mathbf{u}^0,\boldsymbol{\alpha}^0)}$, $i,j=1,...,N_\alpha$ may be considerably smaller than $N_\alpha$, as has been shown in [5-7].

This paper is structured as follows: Section 2 presents the illustrative paradigm benchmark, which models the nonlinear heat conduction in a test section of a within a proposed experimental facility [8,9] for investigating thermal-hydraulics phenomena characterizing the operation and safety of the conceptually-designed G4M reactor [8], a small modular reactor concept cooled by lead-bismuth eutectic (LBE). This paradigm LBE test section benchmark comprises the major ingredients needed for highlighting the salient features involved in applying the *$2^{nd}$-ASAM for nonlinear systems*, yet is sufficiently simple to admit an exact solution, thereby making transparent the mathematical derivations presented in PART I [1]. Section 3 presents the application of the *$2^{nd}$-ASAM* for obtaining the exact expressions of both the first- and second-order sensitivities of the temperature distribution within the test section.

Notably, this application will show that the construction and solution of the *second-level adjoint sensitivity system* (*2$^{nd}$-LASS*) requires very little additional effort beyond the construction of the adjoint sensitivity system needed for computing the first-order sensitivities, and that the actual adjoint computations needed for computing all of the first- and second-order response sensitivities are far less than $N_\alpha$ per response.

In Section 4, the 1$^{st}$- and 2$^{nd}$-order sensitivities are employed to propagate model parameter uncertainties for computing the uncertainties (i.e., variances and skewnesses) in the temperature distribution responses in the heated LBE test section benchmark. Particularizing the general results from [10], Section 4 shows that the 2$^{nd}$-order sensitivities contribute decisively to causing asymmetries in the temperature response distribution. Finally, Section 5 concludes this work by highlighting the most significant results obtained regading the features of the expected temperature distribution in the LBE test section benchmark analyzed herein.

## 2. A PARADIGM NONLINEAR HEAT CONDUCTION PROBLEM

A cylindrical test section for performing heat transfer experiments contains electrically heated rods and is filled with liquid lead-bismuth eutectic (LBE). The length of the cylindrical test section is $\ell = 1.7 m$, and its radius is $a = 15 cm$. The thermal conductivity, denoted as $k(T)$, of the LBE is considered to depend linearly on the temperature, having the functional form

$$k(T) = k_0 (1 + cT), \qquad (1)$$

where the nominal values of the quantities $k_0$ and $c$ are: $k_0^0 = 4.3663 \left[ Wm^{-1} K^{-1} \right]$ and $c^0 = 2.8844 \times 10^{-3} \left[ K^{-1} \right]$. Throughout this work, "nominal values" will be denoted by using the superscript "zero". For simplicity, the electrically heated rods are considered to provide a volumetrically-distributed heat source of nominal strength $Q^0 = 1.11 \times 10^4 \, Wm^{-3}$. The test section is insulated on its lateral surface. The temperature at the bottom of the test section is kept at a constant nominal temperature $T_a^0 = 400 K$. At the top of the section, at $z = \ell/2$, heat is removed by a heat exchanger at a constant heat flux, $q \left[ Wm^{-2} \right]$, having a nominal value $q^0 = 7.44 \times 10^3 \, Wm^{-2}$.

To model mathematically the heat conduction process inside the test section described in the foregoing, it is convenient to take the center of the coordinate systems in the center of the cylinder, so that the test section extends in the axial (vertical) direction from $-\ell/2 \leq z \leq \ell/2$. Since the test section is insulated on its radial surface and since the length of the cylindrical test section is much greater than its radius, the temperature variation in the radial direction can be neglected by comparison to the temperature variations in the axial direction, for the purposes of this illustrative problem. Hence, the axial temperature distribution, $T(z)$, in the LBE can be modeled by the following nonlinear heat conduction model:

$$\frac{d}{dz} \left[ k(T) \frac{dT(z)}{dz} \right] + Q = 0, \quad -\frac{\ell}{2} < z < \frac{\ell}{2}, \qquad (2)$$

$$\left[k(T)\frac{dT}{dz}\right]_{z=\ell/2} = -q, \quad at \ z = \frac{\ell}{2}; \tag{3}$$

$$T(z) = T_a, \quad at \ z = -\frac{\ell}{2}; \tag{4}$$

Notably, *both Eqs. (2) and (3) are nonlinear in $T(z)$, thereby rendering the above heat conduction benchmark problem ideally suited for illustrating the application of the $2^{nd}$-ASAM to the general case of "nonlinear differential equations subject to nonlinear boundary conditions."*

The solution of the above system of nonlinear differential equations can be solved by using Kirchoff's transformation to obtain the solution

$$T(z) = \left[-1 + \sqrt{(1+cT_a)^2 + 2c\tau(z)}\right]\bigg/c, \tag{5}$$

where the function $\tau(z)$ is defined as

$$\tau(z) \equiv \frac{Q}{2k_0}\left(\frac{3\ell^2}{4} + z\ell - z^2\right) - \frac{q}{k_0}\left(z + \frac{\ell}{2}\right). \tag{6}$$

Note that $T(z)$ attains its maximum value, denoted as $T(z_{max})$, at the location

$$z_{max} = \frac{Q\ell - 2q}{2Q}, \tag{7}$$

where it has the expression

$$T(z_{max}) = \left[-1 + \sqrt{(1+cT_a)^2 + 2c\tau(z_{max})}\right]\bigg/c, \tag{8}$$

with $\tau(z_{max})$ given by

$$\tau(z_{max}) = \frac{(Q\ell - q)^2}{2k_0 Q}. \tag{9}$$

For the nominal parameter values provided in the foregoing, the axial variation of the nominal temperature, $T^0(z)$, is depicted in Figure 1. In particular, $T^0(z)$ takes on the following values at the bottom and top, respectively, of the test section: $T^0(-\ell/2) = 400\,K$ and $T^0(\ell/2) = 700\,K$.

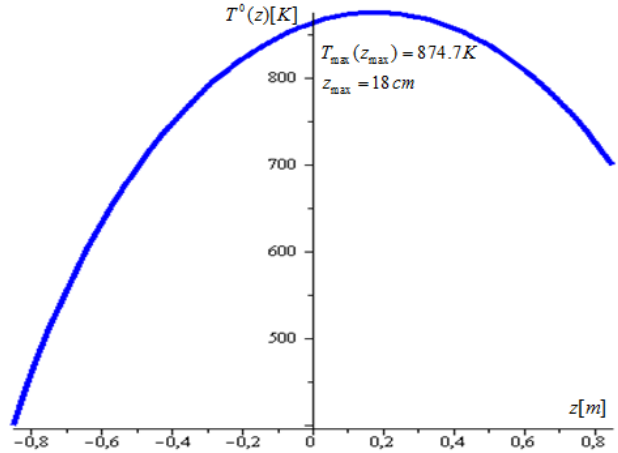

Figure 1: Variation of the nominal temperature, $T(z)$.

It is of interest to install thermocouples within the test section in order to measure the temperatures at various axial locations. A typical thermocouple response would indicate the temperature, $T(z_r)$, at some "response location'" denoted as $z_r$, $-\ell/2 \leq z_r \leq \ell/2$. Such a response is mathematically represented in form

$$R(T;\alpha) \equiv T(z_r) = \int_{-\ell/2}^{\ell/2} T(z)\delta(z - z_r)dz, \tag{10}$$

where $\delta(z - z_r)$ is the customary Dirac delta-functional. A particularly important location for installing a measuring thermocouple is the location $z_{max}$, where $T(z_{max})$ occurs.

Using the notation from Part I [1], the model ("input") parameters in this problem are considered to be the (five) components of the (column) vector

$$\boldsymbol{\alpha} \equiv (Q, q, T_a, k_0, c)^{\dagger}, \tag{11}$$

with nominal values $\boldsymbol{\alpha}^0 \equiv (Q^0, q^0, T_a^0, k_0^0, c^0)^{\dagger}$. The dagger will be used throughout this paper to denote "transposition." The model parameters, $\boldsymbol{\alpha}$, are considered to be afflicted by uncertainties, so they can vary from their nominal values $\boldsymbol{\alpha}^0$ by amounts represented by the component of the "vector of variations", $\mathbf{h}_{\alpha}$, defined as

$$\mathbf{h}_{\alpha} \equiv (\delta Q, \delta q, \delta T_a, \delta k_0, \delta c)^{\dagger}, \tag{12}$$

In practice, the variations $\delta Q, \delta q, \delta T_a, \delta k_0, \delta c$ are usually taken to be the standard deviations quantifying the uncertainties in the respective model parameters.

## 3. APPLICATION OF THE 2$^{nd}$-ASAM FOR COMPUTING THE 1$^{st}$- AND 2$^{nd}$-ORDER SENSITIVITIES

Section 3.1, below, presents the computation of the first-order sensitivities along with their use for computing the first-order contributions to the standard deviation in the temperature distribution, while Section 3.2 presents the computation of the second-order sensitivities along with their use for computing the second-order contributions to the standard deviation and skewness of the temperature distribution.

### *3.1. 2$^{nd}$-ASAM Computation of the First-Order Response Sensitivities*

As shown in PART I [1], the first-order sensitivities of the response $R(\mathbf{e})$, $\mathbf{e} \equiv (T, \boldsymbol{\alpha})$, to the variations $\mathbf{h}_{\alpha}$ are generally obtained by computing the (first-order) G-differential $\delta R(\mathbf{e}^0; \mathbf{h})$ of $R(\mathbf{e})$ at $\mathbf{e}^0 \equiv (T^0, \boldsymbol{\alpha}^0)$, which is defined as

$$\delta R(\mathbf{e}^0;\mathbf{h}) \equiv \frac{d}{d\varepsilon}\{R(\mathbf{e}^0 + \varepsilon\mathbf{h})\}_{\varepsilon=0}, \text{ with } \mathbf{h} \equiv (h_T, \mathbf{h}_\alpha)^\dagger. \tag{13}$$

Applying the above definition to the response defined by Eq. (10) yields the first-order differential, $DR(T^0, \alpha^0; h_\alpha)$, of the response $R(T, \alpha)$:

$$DR(T^0, \boldsymbol{\alpha}^0; h_\alpha) \equiv \frac{d}{d\varepsilon}\left[\int_{-\ell/2}^{\ell/2}[T^0(z) + \varepsilon h_T(z)]\delta(z - z_r)dz\right]_{\varepsilon=0} = \int_{-\ell/2}^{\ell/2} h_T(z)\delta(z - z_r)dz. \tag{14}$$

Next, taking the G-differential of Eqs. (2) - (4) yields

$$\frac{d^2}{dz^2}\left[k(T^0)h_T(z)\right] = -(\delta Q) + (\delta k_0)\frac{Q^0}{k_0^0} - (\delta c)k_0^0\frac{d}{dz}\left[T^0(z)\frac{dT^0(z)}{dz}\right]$$
$$\equiv Q_1(T^0, \boldsymbol{\alpha}; \mathbf{h}_\alpha), \quad z \in (-\ell/2, \ell/2), \tag{15}$$

$$\left\{\frac{d}{dz}[k(T^0)h_T(z)]\right\}_{z=\ell/2} = -(\delta q) + (\delta k_0)\frac{q^0}{k_0^0} - (\delta c)k_0^0\left[T^0(z)\frac{dT^0}{dz}\right]_{z=\ell/2}$$
$$\equiv q_1(T^0, \boldsymbol{\alpha}; \mathbf{h}_\alpha), \text{ at } z = \ell/2; \tag{16}$$

$$h_T(z)\big|_{z=-\ell/2} = (\delta T_a), \text{ at } z = -\ell/2. \tag{17}$$

Equations (15)-(17) correspond to Eqs. (13)-(14) in [1]. Define the Hilbert space $\mathcal{H}(\Omega) \equiv \mathcal{L}_2(\Omega)$ to consist of all square integrable functions $f(z)$ defined on the domain $z \in \Omega \equiv [-\ell/2, \ell/2]$ and endowed with the inner product [for two functions, $f_1(z)$ and $f_1(z)$] as follows:

$$\langle f_1(z), f_2(z)\rangle \equiv \int_{-\ell/2}^{\ell/2} f_1(z)f_2(z)dz. \tag{18}$$

Recall from the general theory presented in [1] that the may require two distinct Hilber spaces, denoted in [1] as $\mathcal{H}_u(\Omega_x)$ and $\mathcal{H}_Q(\Omega_x)$, respectively. For this simple iilustrative problem, however, both of these Hilbert spaces coincide with $\mathcal{L}_2(\Omega)$.

Forming now the inner product of Eq. (15) with a yet undefined function $\Psi(z) \in \mathcal{L}_2(\Omega)$ and integrating the resulting equations twice by parts to transfer the differential operations from $h_T(z)$ to $\Psi(z)$ yields:

$$\int_{-\ell/2}^{\ell/2} \Psi(z) \frac{d^2}{dz^2}\left[k(T^0)h_T(z)\right] dz = \int_{-\ell/2}^{\ell/2} \Psi(z) Q_1(T^0, \boldsymbol{\alpha}^0; \mathbf{h}_\alpha) dz$$
$$= \left\{\Psi(z) \frac{d}{dz}\left[k(T^0)h_T(z)\right] - h_T(z) k(T^0) \frac{d\Psi}{dz}\right\}_{z=-\ell/2}^{z=\ell/2} + \int_{-\ell/2}^{\ell/2} h_T(z)\left[k(T^0) \frac{d^2\Psi}{dz^2}\right] dz. \quad (19)$$

The above Eq. (19) corresponds to Eq. (15) in [1]. Applying the principles outlined in Part I to Eq. (19), yields the following *first-level adjoint sensitivity system (1st-LASS)*:

$$k(T^0) \frac{d^2\Psi(z)}{dz^2} = \delta(z - z_r), \quad z \in (-\ell/2, \ell/2), \quad (20)$$

$$\left.\frac{d\Psi}{dz}\right|_{z=\ell/2} = 0, \quad at\ z = \ell/2, \quad (21)$$

$$\Psi(z) = 0, \quad at\ z = -\ell/2. \quad (22)$$

The *1st-LASS* above can be readily solved to obtain the *1st- level adjoint function*

$$\Psi(z) = \frac{(z - z_r) H(z - z_r) - z - \ell/2}{k\left[T^0(z_r)\right]}. \quad (23)$$

where $H(z)$ is the customary Heaviside unit-step functional, defined as $H(z) = 0\ if\ z < 0$ and $H(z) = 1\ if\ z \geq 0$. Note that the *1st- level adjoint function* $\Psi(z)$ can also be interpreted to be the Green`s function

$$G(z, z_r) = \frac{1}{k[T(z_r)]} \left[ (z - z_r) H(z - z_r) + z - \ell/2 \right], \tag{24}$$

for the *1$^{st}$-LASS*, i.e., Eqs. (20) - (22), since the point $z_r \in [-\ell/2, \ell/2]$ is arbitrary.

As shown in the general theory in PART I [1], the *1$^{st}$-LASS* depends on the response, cf. Eq.(14), which provides the source term $\delta(z - z_r)$, as shown in Eq. (20). Note that this source term does not belong to the Hilbert space $\mathscr{L}_2(\Omega)$ but belongs to the Sobolev space $\mathscr{H}_0^1(\Omega) \subset \mathscr{L}_2(\Omega)$, i.e., $\Sigma_d^0 \delta(x - b) \in \mathscr{H}_0^1 \subset \mathscr{L}_2(\Omega)$, as usually encountered when computing Green's functions. By using the well known Lax-Milgram Lemma, it can be shown that the bilinear form on the right side of the last equality in Eq. (19) coercive, so that the *1$^{st}$-LASS* can be solved uniquely, as has been done to obtain the expression for the adjoint function $\Psi(z) \in \mathscr{H}(\Omega) = \mathscr{H}_0^1(\Omega) \subset \mathscr{L}_2(\Omega)$ shown in Eq. (23). The mathematical technicalities requiring the use of Sobolev spaces stemming from the consideration of distributions, as encountered above for the *first-level adjoint sensitivity system* (*1$^{st}$-LASS*), will also arise in the construction of the *second-level adjoint sensitivity systems* (*2$^{nd}$-LASS*), as will be seen in the remainder of this work. Thus, even though the foregoing mathematical technicalities will not be repeated in the sequel, all of the solutions to such *2$^{nd}$-LASS* should be interpreted in the "weak sense," in the appropriate Sobolev space. No confusion should arise, however, since the respective solutions for the *2$^{nd}$-LASS* will be unique and will be obtained explicitly, just as it was in Eq. (23) for the first-level adjoint function $\Psi(z)$.

Using the results in Eqs. (15) - (22) in Eq. (14) transforms the latter into the form:

$$DR(T^0, \boldsymbol{\alpha}^0; \mathbf{h}_\alpha) \equiv \int_{-\ell/2}^{\ell/2} \Psi(z) Q_1(T^0, \boldsymbol{\alpha}^0; \mathbf{h}_\alpha) dz \\ -(\delta T_a)\left[ k(T^0) \frac{d\Psi}{dz} \right]_{z=-\ell/2} - \Psi(z = \ell/2) q_1(T^0, \boldsymbol{\alpha}^0; \mathbf{h}_\alpha). \tag{25}$$

Using Eqs. (15) and (16), respectively, to replace the expressions of $Q_1(T^0, \boldsymbol{\alpha}^0; \mathbf{h}_\alpha)$ and $q_1(T^0, \boldsymbol{\alpha}^0; \mathbf{h}_\alpha)$ in Eq. (35) yields

$$DR(T^0, \boldsymbol{\alpha}^0; \mathbf{h}_\alpha) = \frac{\partial T(z_r)}{\partial Q}(\delta Q) + \frac{\partial T(z_r)}{\partial q}(\delta q) + \frac{\partial T(z_r)}{\partial T_a}(\delta T_a) + \frac{\partial T(z_r)}{\partial k_0}(\delta k_0) + \frac{\partial T(z_r)}{\partial c}(\delta c) \tag{26}$$

where the the 1$^{st}$-order sensitivities of $T(z_r)$ to the repective model parameters have the following expressions:

$$S_1(T,\Psi;\boldsymbol{\alpha}) \equiv \frac{\partial T(z_r)}{\partial Q} = -\int_{-\ell/2}^{\ell/2} \Psi(z)dz, \tag{27}$$

$$S_2(T,\Psi;\boldsymbol{\alpha}) \equiv \frac{\partial T(z_r)}{\partial q} = \Psi(\ell/2), \tag{28}$$

$$S_3(T,\Psi;\boldsymbol{\alpha}) \equiv \frac{\partial T(z_r)}{\partial T_a} = -\left[k(T^0)\frac{d\Psi}{dz}\right]_{z=-\ell/2}, \tag{29}$$

$$S_4(T,\Psi;\boldsymbol{\alpha}) \equiv \frac{\partial T(z_r)}{\partial k_0} = \frac{1}{k_0^0}\left[Q^0\int_{-\ell/2}^{\ell/2}\Psi(z)dz - q^0\Psi(\ell/2)\right], \tag{30}$$

$$\begin{aligned}S_5(T,\Psi;\boldsymbol{\alpha}) \equiv \frac{\partial T(z_r)}{\partial c} &= -\int_{-\ell/2}^{\ell/2}\Psi(z)k_0^0\frac{d}{dz}\left[T^0(z)\frac{dT^0(z)}{dz}\right]dz \\ &+ \Psi(\ell/2)k_0^0\left[T^0(z)\frac{dT^0}{dz}\right]_{z=\ell/2} = k_0^0\frac{T_a^2-[T^0(z_r)]^2}{2k[T^0(z_r)]}.\end{aligned} \tag{31}$$

The above expressions indicate that all of the 1$^{st}$-order sensitivities can be computed exacly and efficiently, using quadratures (integrations), once the adjoint function has been obtained by solving the *1$^{st}$-LASS*. Thus, replacing the expression of $\Psi(z)$ given in Eq. (23) in Eqs. (27) – (31) and carrying out the respective integrations over $z$ yields the following evaluated expressions for the 1$^{st}$-order sensitivities of the response $T(z_r)$ to the model parameters:

$$\frac{\partial T(z_r)}{\partial Q} = \frac{1}{2k[T^0(z_r)]}\left(\frac{3\ell^2}{4} + z_r\ell - z_r^2\right), \tag{32}$$

$$\frac{\partial T(z_r)}{\partial q} = -\frac{z_r + \ell/2}{k[T^0(z_r)]}, \tag{33}$$

$$\frac{\partial T(z_r)}{\partial T_a} = \frac{1+c^0 T_a^0}{1+c^0 T^0(z_r)}, \tag{34}$$

$$\frac{\partial T(z_r)}{\partial k_0} = \frac{1}{k_0^0 k\left[T^0(z_r)\right]} \left[ -\frac{Q^0}{2}\left(\frac{3\ell^2}{4}+z_r\ell-z_r^2\right) + q^0\left(z_r+\frac{\ell}{2}\right)\right] = -\frac{\tau(z_r)}{k\left[T^0(z_r)\right]}, \tag{35}$$

$$\frac{\partial T(z_r)}{\partial c} = k_0^0 \frac{T_a^2-\left[T^0(z_r)\right]^2}{2k\left[T^0(z_r)\right]} = \frac{T^0(z_r)-\tau(z_r)-T_a^0}{c^0\left[1+c^0 T^0(z_r)\right]}. \tag{36}$$

One of the main uses of sensitivities is for ranking the relative importance of parameter variations in influencing variations in responses. Relative sensitivities are used for this purpose, since they are dimensionless numbers. The relative sensitivity of a response $R(\mathbf{e})$ to the $i^{th}$-parameter, $\alpha_i$, is defined as $S_i^{rel} \triangleq (\partial R/\partial \alpha_i)_{\mathbf{e}^0}\left[\alpha_i^0/R(\mathbf{e}^0)\right]$. The relative sensitivities of $T(z_r)$ are depicted in Figures 2-6 as functions of the arbitrary location $z_r \in [-\ell/2, \ell/2]$.

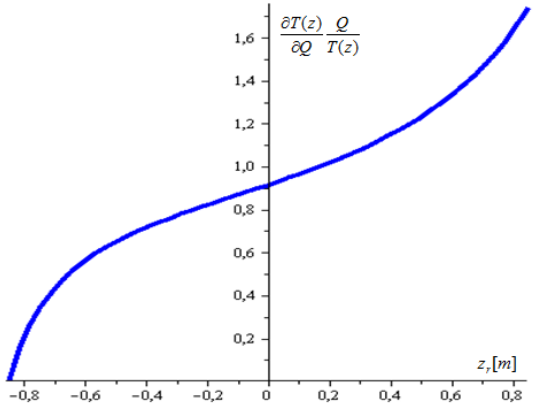 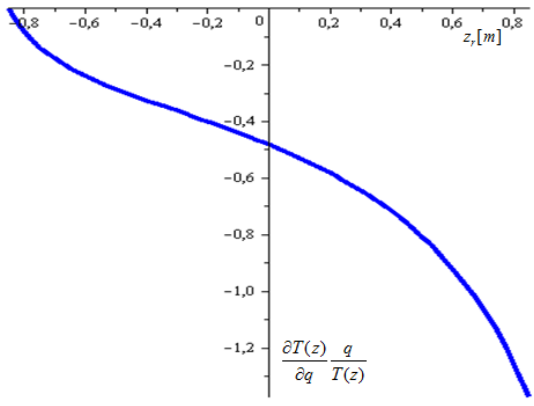

Figure 2: Relative sensitivity of $T(z_r)$ to $Q$.    Figure 3: Relative sensitivity of $T(z_r)$ to $q$.

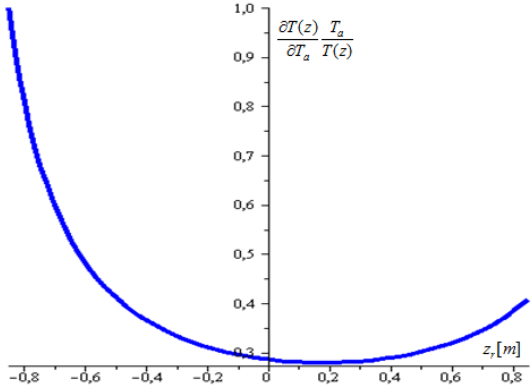
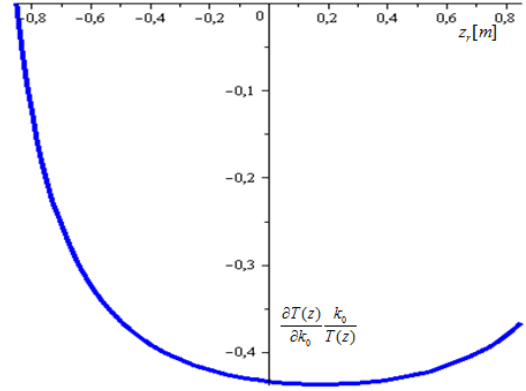

Figure 4: Relative sensitivity of $T(z_r)$ to $T_a$.   Figure 5: Relative sensitivity of $T(z_r)$ to $k_0$.

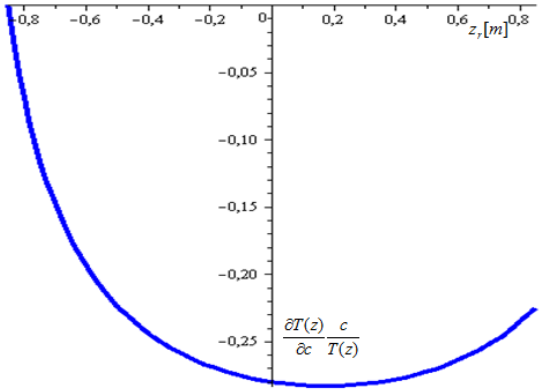

Figure 6: Relative sensitivity of $T(z_r)$ to $c$.

As Figures 2 through 6 indicate, the relative importance of the various model parameters depends on the position $z_r$. Ranking them by the largest attainable maximum absolute values of their relative sensitivities, the importance of the model parameters is as follows: the heat source $Q$; the boundary heat flux $q$; the ambient temperature $T_a$; the linear heat conductivity coefficient $k_0$; and the nonlinear heat conductivity coefficient $c$.

### 3.2. 2$^{nd}$-ASAM Computation of the Second-Order Response Sensitivities

As discussed in the general theory presented in PART I [1], the fundamental philosophical consideration and starting point for computing the second-order response sensitivities, $S_{ij}(\alpha^0) \triangleq \partial^2 R / \partial \alpha_i \partial \alpha_j = \partial^2 R / \partial \alpha_j \partial \alpha_i$, was to consider the first-order sensitivities to be responses of the form $S_i(T, \Psi; \alpha)$, $i = 1, 2, 3, 4, 5$. This fact was explicitly indicated in the respective definitions, provided by Eqs. (27) through (31). Based on this fundamental consideration, the $2^{nd}$-ASAM proceeds by computing the first-order G-differential, $\delta S_i(T, \Psi; \alpha)$, of each of the functionals $S_i(T, \Psi; \alpha)$ at the point $(T^o, \Psi^0; \alpha^0)$ of nominal values, using the definition of the G-differential, namely:

$$\delta S_i(T^o, \Psi^0; \alpha^0; h_T, h_\Psi, \mathbf{h}_\alpha) \equiv \left\{ \frac{d}{d\varepsilon} \left[ S_i(T^o + \varepsilon h_T, \Psi^0 + \varepsilon h_\Psi, \alpha^0 + \varepsilon \mathbf{h}_\alpha) \right] \right\}_{\varepsilon=0} \tag{37}$$

for an arbitrary scalar $\varepsilon \in \mathcal{F}$, and vectors $(h_T, h_\Psi, \mathbf{h}_\alpha) \in \mathcal{H}(\Omega) \times \mathcal{H}(\Omega) \times \mathcal{H}_\alpha$. For our illustrative example, there will be $N_\alpha (N_\alpha + 1)/2 = 15$ distinct second-order derivatives, due to the symmetry property $\partial^2 R / \partial \alpha_i \partial \alpha_j = \partial^2 R / \partial \alpha_j \partial \alpha_i$. The sensitivities $S_{ij}(\alpha^0)$ will be computed next, in the order of parameter importance/ranking as discussed in the previous sub-section.

*3.2.1. Computation of the Second-Order Response Sensitivities $S_{1i} \triangleq \partial^2 T(z_r)/(\partial Q \partial \alpha_i)$*

Applying the definition shown in Eq. (37) to Eq. (27) yields the G-differential, $DS_1(T^o, \Psi^0; \alpha^0; h_T, h_\Psi, \mathbf{h}_\alpha)$, of the first-order sensitivity $S_1(T^o, \Psi^0; \alpha^0)$, in the form

$$DS_1(T^o, \Psi^0; \alpha^0; h_T, h_\Psi, \mathbf{h}_\alpha) = \frac{d}{d\varepsilon} \left\{ -\int_{-\ell/2}^{\ell/2} \Psi(z) dz \right\}_{\varepsilon=0} = -\int_{-\ell/2}^{\ell/2} h_\Psi(z) dz . \tag{38}$$

The function $h_\Psi$ in Eq. (38) is the solution of the system of equations obtained by G-differentiating Eqs. (20) - (22), namely

$$k\left(T^{0}\right)\frac{d^{2}h_{\Psi}(z)}{dz^{2}}+k_{0}^{0}c^{0}h_{T}(z)\frac{d^{2}\Psi^{0}(z)}{dz^{2}}$$
$$=-\left\{(\delta k_{0})\left[1+c^{0}T^{0}(z)\right]+(\delta c)k_{0}^{0}T^{0}(z)\right\}\frac{d^{2}\Psi^{0}(z)}{dz^{2}},\ z\in(-\ell/2,\ \ell/2), \tag{39}$$

$$\left.\frac{dh_{\Psi}}{dz}\right|_{z=\ell/2}=0,\ at\ z=\ell/2, \tag{40}$$

$$h_{\Psi}(z)=0,\ at\ z=-\ell/2. \tag{41}$$

The quantity $d^{2}\Psi^{0}(z)/dz^{2}$ in Eq. (39) can be replaced by using Eq. (20) to obtain the following form of Eq. (39)

$$k\left(T^{0}\right)\frac{d^{2}h_{\Psi}(z)}{dz^{2}}+h_{T}(z)\frac{c^{0}\delta(z-z_{r})}{1+c^{0}T^{0}(z)}$$
$$=-\left[\frac{(\delta k_{0})}{k_{0}^{0}}+(\delta c)\frac{T^{0}(z)}{1+c^{0}T^{0}(z)}\right]\delta(z-z_{r})\equiv Q_{12}\left(T^{0},\boldsymbol{\alpha}^{0};\mathbf{h}_{\alpha}\right),\ z\in(-\ell/2,\ \ell/2). \tag{42}$$

As Eq. (38) indicates, the entire contribution to $DS_{1}\left(T^{o},\Psi^{0};\boldsymbol{\alpha}^{0};h_{T},h_{\Psi},\mathbf{h}_{\alpha}\right)$ comes from the "indirect-effect" term; there is no "direct-effect" term contribution to $DS_{1}\left(T^{o},\Psi^{0};\boldsymbol{\alpha}^{0};h_{T},h_{\Psi},\mathbf{h}_{\alpha}\right)$. Applying next the general theoretical considerations leading to Eq. (34) of PART I [1], Eqs. (15) and (42) are written in the following block-matrix form:

$$\begin{pmatrix}\dfrac{d^{2}}{dz^{2}}\left\{k\left(T^{0}\right)[\ ]\right\} & 0 \\ \dfrac{c^{0}\delta(z-z_{r})}{1+c^{0}T^{0}(z)} & k\left(T^{0}\right)\dfrac{d^{2}[\ ]}{dz^{2}}\end{pmatrix}\begin{pmatrix}h_{T}(z)\\ h_{\Psi}(z)\end{pmatrix}=\begin{pmatrix}Q_{1}\left(T^{0},\boldsymbol{\alpha};\mathbf{h}_{\alpha}\right)\\ Q_{12}\left(T^{0},\boldsymbol{\alpha}^{0};\mathbf{h}_{\alpha}\right)\end{pmatrix},\ z\in(-\ell/2,\ \ell/2). \tag{43}$$

Following the procedure outlined in [1], introduce the vector $\boldsymbol{\Psi}_{1}^{(2)}\equiv\left(\Psi_{11}^{(2)},\ \Psi_{12}^{(2)}\right)\in\mathscr{H}(\Omega)\times\mathscr{H}(\Omega)$ and define the inner product between two functions $\mathbf{h}^{(2)}(z)\equiv\left[h_{T}(z),\ h_{\Psi}(z)\right]$ and $\boldsymbol{\Psi}_{1}^{(2)}\equiv\left[\Psi_{11}^{(2)}(z),\ \Psi_{12}^{(2)}(z)\right]$ as follows:

$$\left\langle \mathbf{h}^{(2)}(z), \mathbf{\Psi}_1^{(2)} \right\rangle \equiv \int_{-\ell/2}^{\ell/2} \left[ h_T(z) \Psi_{11}^{(2)}(z) + h_\Psi(z) \Psi_{12}^{(2)}(z) \right] dz. \tag{44}$$

Following the sequence of operations leading to Eq. (37) of [1], form the inner product of $\mathbf{\Psi}_1^{(2)} \equiv \left[ \Psi_{11}^{(2)}(z), \Psi_{12}^{(2)}(z) \right]$ with Eq. (43) to obtain the following sequence of equalities:

$$\begin{aligned}
&\left\langle \left[ \Psi_{11}^{(2)}(z), \Psi_{12}^{(2)}(z) \right] \begin{pmatrix} \dfrac{d^2}{dz^2}\{k(T^0)[\ ]\} & 0 \\ \dfrac{c^0 \delta(z-z_r)}{1+c^0 T^0(z)} & k(T^0)\dfrac{d^2[\ ]}{dz^2} \end{pmatrix} \begin{pmatrix} h_T(z) \\ h_\Psi(z) \end{pmatrix} \right\rangle \\
&= \left\langle \left[ \Psi_{11}^{(2)}(z), \Psi_{12}^{(2)}(z) \right] \begin{pmatrix} Q_1(T^0, \boldsymbol{\alpha}; \mathbf{h}_\alpha) \\ Q_{12}(T^0, \boldsymbol{\alpha}^0; \mathbf{h}_\alpha) \end{pmatrix} \right\rangle \\
&= \left\langle \left[ h_T(z), h_\Psi(z) \right] \begin{pmatrix} k(T^0)\dfrac{d^2[\ ]}{dz^2} & \dfrac{c^0 \delta(z-z_r)}{1+c^0 T^0(z)} \\ 0 & \dfrac{d^2}{dz^2}\{k(T^0)[\ ]\} \end{pmatrix} \begin{pmatrix} \Psi_{11}^{(2)}(z) \\ \Psi_{12}^{(2)}(z) \end{pmatrix} \right\rangle \\
&\quad + \left\{ P_2(T^0, \boldsymbol{\alpha}^0; h_T, h_\Psi, \mathbf{h}_\alpha; \Psi_{11}^{(2)}, \Psi_{12}^{(2)}) \right\}_{\partial \Omega},
\end{aligned} \tag{45}$$

where the bilinear concomitant $\left\{ P_2(T^0, \boldsymbol{\alpha}^0; h_T, h_\Psi, \mathbf{h}_\alpha; \Psi_{11}^{(2)}, \Psi_{12}^{(2)}) \right\}_{\partial \Omega}$ has the form

$$\left\{ P_2(T^0, \boldsymbol{\alpha}^0; h_T, h_\Psi; \Psi_{11}^{(2)}, \Psi_{12}^{(2)}) \right\}_{\partial \Omega} \equiv \left\{ \Psi_{11}^{(2)}(z) \frac{d\left[ k(T^0) h_T(z) \right]}{dz} \right.$$

$$\left. -k(T^0) h_T(z) \frac{d\Psi_{11}^{(2)}(z)}{dz} + \Psi_{12}^{(2)}(z) k(T^0) \frac{dh_\Psi(z)}{dz} - h_\Psi(z) \frac{d\left[ k(T^0) \Psi_{12}^{(2)}(z) \right]}{dz} \right\}\Bigg|_{z=-\ell/2}^{z=\ell/2}. \tag{46}$$

As generally shown in [1], the definition of the 2$^{nd}$-level adjoint function $\mathbf{\Psi}_1^{(2)} \equiv \left( \Psi_{11}^{(2)}, \Psi_{12}^{(2)} \right)$ is now completed by requiring the term on the right side on the last equality in Eq. (45) to represent the same functional as the right side of Eq. (38), which yields the following "*second-level adjoint sensitivity system* (2$^{nd}$-*LASS*):"

$$\begin{pmatrix} k(T^0)\dfrac{d^2[\ ]}{dz^2} & \dfrac{c^0\delta(z-z_r)}{1+c^0T^0(z)} \\ 0 & \dfrac{d^2}{dz^2}\{k(T^0)[\ ]\} \end{pmatrix} \begin{pmatrix} \Psi_{11}^{(2)}(z) \\ \Psi_{12}^{(2)}(z) \end{pmatrix} = \begin{pmatrix} 0 \\ -1 \end{pmatrix} \qquad (47)$$

As generally discussed in [1], the boundary conditions for the above $2^{nd}$-LASS are obtained by using the boundary conditions given in Eqs. (16), (17), (40) and (41) in Eq. (46) to eliminate all unknown values of $h_T(z)$ and $h_\Psi(z)$, respectively, in the bilinear concomitant $\{P_2(T^0,\mathbf{\alpha}^0;h_T,h_\Psi;\Psi_{11}^{(2)},\Psi_{12}^{(2)})\}_{\partial\Omega}$. These considerations lead to the following boundary conditions for the $2^{nd}$-LASS:

$$\left.\dfrac{d\Psi_{11}^{(2)}}{dz}\right|_{z=\ell/2} = 0, \quad at\ z=\ell/2, \qquad (48)$$

$$\Psi_{11}^{(2)}(z) = 0, \quad at\ z = -\ell/2. \qquad (49)$$

$$\left.\dfrac{d\left[k(T^0)\Psi_{12}^{(2)}\right]}{dz}\right|_{z=\ell/2} = 0, \quad at\ z = \ell/2, \qquad (50)$$

$$\Psi_{12}^{(2)}(z) = 0, \quad at\ z = -\ell/2. \qquad (51)$$

Inserting the above boundary conditions together with those those given in Eqs. (16), (17), (40) and (41) into Eq. (46) reduces the bilinear concomitant $\{P_2(T^0,\mathbf{\alpha}^0;h_T,h_\Psi;\Psi_{11}^{(2)},\Psi_{12}^{(2)})\}_{\partial\Omega}$ to the quantity $\hat{P}_2(T^0,\mathbf{\alpha}^0;\Psi_{11}^{(2)},\Psi_{12}^{(2)};\mathbf{h}_\alpha)$, which has the following form:

$$\hat{P}_2(T^0,\mathbf{\alpha}^0;\Psi_{11}^{(2)},\Psi_{12}^{(2)};\mathbf{h}_\alpha) \equiv q_1(T^0,\mathbf{\alpha};\mathbf{h}_\alpha)\left[\Psi_{11}^{(2)}(z)\right]_{z=\ell/2} + (\delta T_a)\left[k(T^0)\dfrac{d\Psi_{11}^{(2)}(z)}{dz}\right]_{z=-\ell/2}. \qquad (52)$$

Finally, using Eqs. (52), (47) and (38) in Eq. (45) leads to the following expression for the second-order differential expression $DS_1(T^0;\mathbf{\alpha}^0;h_T,h_\Psi,\mathbf{h}_\alpha)$

$$DS_1\left(T^0;\boldsymbol{\alpha}^0;h_T,h_\Psi,\mathbf{h}_\alpha\right) = \int_{-\ell/2}^{\ell/2} \left[\Psi_{11}^{(2)}(z)Q_1\left(T^0,\boldsymbol{\alpha}^0;\mathbf{h}_\alpha\right) + \Psi_{12}^{(2)}(z)Q_{12}\left(T^0,\boldsymbol{\alpha}^0;\mathbf{h}_\alpha\right)\right]dz$$
$$-q_1\left(T^0,\boldsymbol{\alpha};\mathbf{h}_\alpha\right)\left[\Psi_{11}^{(2)}(z)\right]_{z=\ell/2} - \left(\delta T_a\right)\left[k\left(T^0\right)\frac{d\Psi_{11}^{(2)}(z)}{dz}\right]_{z=-\ell/2}. \tag{53}$$

Solving Eqs. (47), (50) and (51) yields the following expression for the 2$^{\text{nd}}$-level adjoint function $\Psi_{12}^{(2)}(z)$:

$$\Psi_{12}^{(2)}(z) = \frac{1}{2k\left[T^0(z)\right]}\left(\frac{3\ell^2}{4} + z\ell - z^2\right). \tag{54}$$

Solving Eqs. (47) - (49) yields the following expression for the 2$^{\text{nd}}$-level adjoint function $\Psi_{11}^{(2)}(z)$:

$$\Psi_{11}^{(2)}(z) = C\left(T^0, z_r;\boldsymbol{\alpha}^0\right)\frac{1}{2}\left(\frac{3\ell^2}{4} + z_r\ell - z_r^2\right)\left[z + \frac{\ell}{2} - (z - z_r)H(z - z_r)\right]. \tag{55}$$

with

$$C\left(T^0, z_r;\boldsymbol{\alpha}^0\right) \equiv \frac{k_0^0 c^0}{k^3\left[T^0(z_r)\right]}. \tag{56}$$

Replacing $Q_1\left(T^0,\boldsymbol{\alpha}^0;\mathbf{h}_\alpha\right)$ and $Q_{12}\left(T^0,\boldsymbol{\alpha}^0;\mathbf{h}_\alpha\right)$ with the corresponding expressions from Eqs. (15) and (42), respectively, yields the following expression for the differential $DS_1\left(T^0;\boldsymbol{\alpha}^0;h_T,h_\Psi,\mathbf{h}_\alpha\right)$:

$$DS_1\left(T^0;\boldsymbol{\alpha}^0;h_T,h_\Psi,\mathbf{h}_\alpha\right) = \frac{\partial^2 T(z_r)}{\partial Q^2}(\delta Q) + \frac{\partial^2 T(z_r)}{\partial Q \partial q}(\delta q)$$
$$+ \frac{\partial^2 T(z_r)}{\partial Q \partial T_a}(\delta T_a) + \frac{\partial^2 T(z_r)}{\partial Q \partial k_0}(\delta k_0) + \frac{\partial^2 T(z_r)}{\partial Q \partial c}(\delta c), \tag{57}$$

where

$$\frac{\partial^2 T(z_r)}{\partial Q^2} = -\int_{-\ell/2}^{\ell/2} \Psi_{11}^{(2)}(z)dz, \tag{58}$$

$$\frac{\partial^2 T(z_r)}{\partial Q \partial q} = \left[\Psi_{11}^{(2)}(z)\right]_{z=\ell/2}, \tag{59}$$

$$\frac{\partial^2 T(z_r)}{\partial Q \partial T_a} = -\left[k(T^0)\frac{d\Psi_{11}^{(2)}(z)}{dz}\right]_{z=-\ell/2}, \tag{60}$$

$$\frac{\partial^2 T(z_r)}{\partial Q \partial k_0} = \frac{1}{k_0^0}\left[Q^0 \int_{-\ell/2}^{\ell/2} \Psi_{11}^{(2)}(z)dz - \int_{-\ell/2}^{\ell/2} \Psi_{12}^{(2)}(z)\delta(z-z_r)dz - q^0\left[\Psi_{11}^{(2)}(z)\right]_{z=\ell/2}\right], \tag{61}$$

$$\begin{aligned}\frac{\partial^2 T(z_r)}{\partial Q \partial c} &= -k_0^0 \int_{-\ell/2}^{\ell/2} \Psi_{11}^{(2)}(z)\frac{d}{dz}\left[T^0(z)\frac{dT^0(z)}{dz}\right]dz \\ &\quad - \int_{-\ell/2}^{\ell/2} \Psi_{12}^{(2)}(z)\frac{T^0(z)}{1+c^0 T^0(z)}\delta(z-z_r)dz + k_0^0\left[\Psi_{11}^{(2)}(z)T^0(z)\frac{dT^0}{dz}\right]_{z=\ell/2} \\ &= -k_0^0\left\{\frac{T_a^2}{2}\left[\frac{d\Psi_{11}^{(2)}}{dz}\right]_{z=-\ell/2} + \frac{1}{2}\int_{-\ell/2}^{\ell/2}[T^0(z)]^2\frac{d^2\Psi_{11}^{(2)}}{dz^2}dz + \frac{\Psi_{12}^{(2)}(z_r)T^0(z_r)}{k[T^0(z_r)]}\right\}.\end{aligned} \tag{62}$$

Replacing Eqs. (54) and (55) in the above expressions and carrying out the integrations over $z$ yields the following expressions for the above 2$^{nd}$-order sensitivities:

$$\frac{\partial^2 T(z_r)}{\partial Q^2} = -\frac{1}{4}C(T^0, z_r; \boldsymbol{\alpha}^0)\left(\frac{3\ell^2}{4} + z_r\ell - z_r^2\right)^2, \tag{63}$$

$$\frac{\partial^2 T(z_r)}{\partial Q \partial q} = \frac{1}{2}C(T^0, z_r; \boldsymbol{\alpha}^0)\left(\frac{3\ell^2}{4} + z_r\ell - z_r^2\right)\left(z_r + \frac{\ell}{2}\right), \tag{64}$$

$$\frac{\partial^2 T(z_r)}{\partial Q \partial T_a} = -C(T^0, z_r; \boldsymbol{\alpha}^0)\frac{k(T_a^0)}{2}\left(\frac{3\ell^2}{4} + z_r\ell - z_r^2\right), \tag{65}$$

$$\frac{\partial^2 T(z_r)}{\partial Q \partial k_0} = C(T^0, z_r; \boldsymbol{\alpha}^0)\frac{1}{2}\left\{\tau^0(z_r) - \frac{1}{c^0}\left[1 + c^0 T^0(z_r)\right]^2\right\}\left(\frac{3\ell^2}{4} + z_r\ell - z_r^2\right), \tag{66}$$

$$\frac{\partial^2 T(z_r)}{\partial Q \partial c} = C(T^0, z_r; \mathbf{\alpha}^0) \frac{k_0^0}{2c^0} \left\{ T_a^0 + \tau^0(z_r) - 2T^0(z_r) - c^0 \left[ T^0(z_r) \right]^2 \right\}$$
$$\times \left( \frac{3\ell^2}{4} + z_r \ell - z_r^2 \right). \quad (67)$$

The relative sensitivities The 2$^{nd}$-order relative sensitivities $(S_{1i})_{rel} \triangleq \frac{\partial^2 T(z_r)}{\partial Q \partial \alpha_i} \frac{Q \alpha_i}{T(z_r)}$ are depicted in Figures 7 - 11 as functions of the arbitrary location $z_r \in [-\ell/2, \ell/2]$.

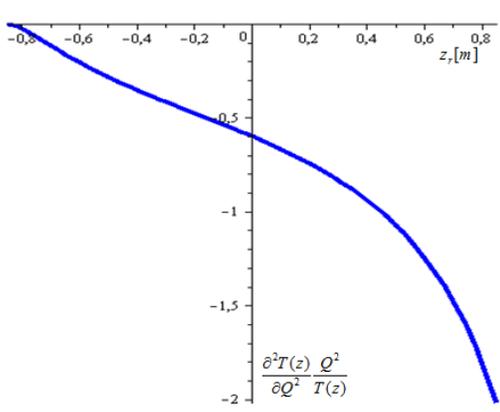 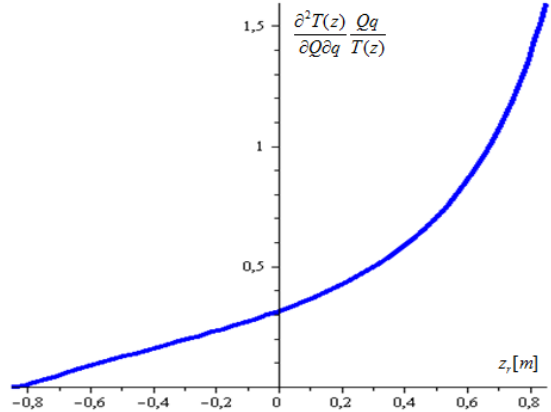

Figure 7: Relative sensitivity $\frac{\partial^2 T(z_r)}{\partial Q^2} \frac{Q^2}{T(z_r)}$. Figure 8: Relative sensitivity $\frac{\partial^2 T(z_r)}{\partial Q \partial q} \frac{Qq}{T(z_r)}$

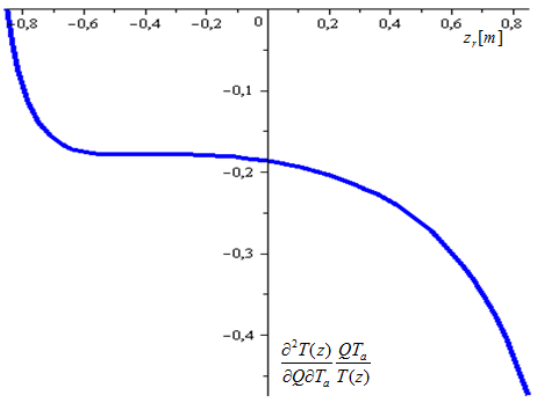 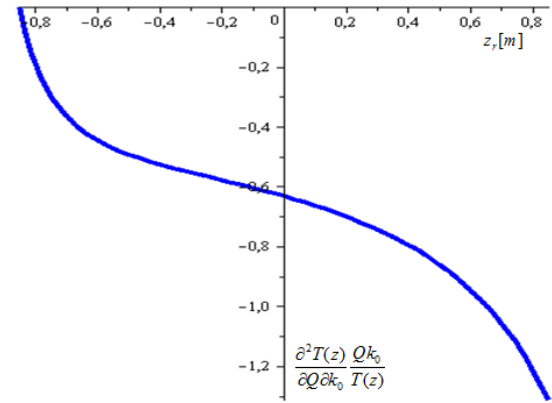

Figure 9: Relative sensitivity $\frac{\partial^2 T(z_r)}{\partial Q \partial T_a} \frac{Q T_a}{T(z_r)}$. Figure 10: Relative sensitivity $\frac{\partial^2 T(z_r)}{\partial Q \partial k_0} \frac{Q k_0}{T(z_r)}$

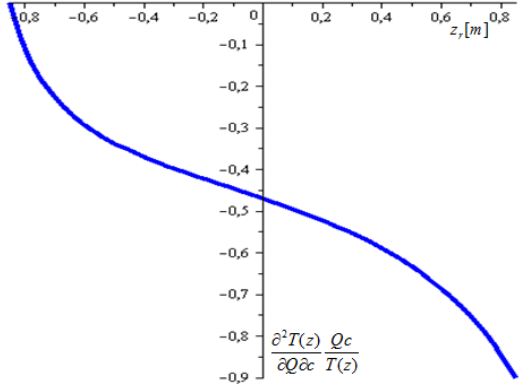

Figure 11: Relative sensitivity $\dfrac{\partial^2 T(z_r)}{\partial Q \partial c}\dfrac{Qc}{T(z_r)}$

As Figures 7 through 11 indicate, the the relative importance (magnitude) of the 2$^{nd}$-order relative sensitivities $(S_{1i})_{rel} \triangleq \dfrac{\partial^2 T(z_r)}{\partial Q \partial \alpha_i}\dfrac{Q\alpha_i}{T(z_r)}$ depends on the position $z_r$. Ranking them by the largest attainable maximum absolute values of these relative sensitivities, the importance of the model parameters is as follows: the heat source $Q$; the boundary heat flux $q$; the linear heat conductivity coefficient $k_0$; the ambient temperature $T_a$; and the nonlinear heat conductivity coefficient $c$. By comparison to the importance ranking of the 1$^{st}$-order relative sensitivities, the model parameters $k_0$ and $T_a$ have "switched places," in that the relative sensitivity $\dfrac{\partial^2 T(z_r)}{\partial Q \partial k_0}\dfrac{Qk_0}{T(z_r)}$ has become about 3 times as large as the relative sensitivity $\dfrac{\partial^2 T(z_r)}{\partial Q \partial T_a}\dfrac{QT_a}{T(z_r)}$. Furthermore, the relative sensitivity $\dfrac{\partial^2 T(z_r)}{\partial Q \partial c}\dfrac{Qc}{T(z_r)}$ has become third in the importance ranking of the partial sensitivities $(S_{1i})_{rel} \triangleq \dfrac{\partial^2 T(z_r)}{\partial Q \partial \alpha_i}\dfrac{Q\alpha_i}{T(z_r)}$.

*3.2.2. Computation of the Second-Order Response Sensitivities* $S_{2i} \triangleq \partial^2 T(z_r)/(\partial q \partial \alpha_i)$

Applying the definition shown in Eq. (37) to Eq. (28) yields the G-differential, $DS_2(T^o, \Psi^0; \boldsymbol{\alpha}^0; h_T, h_\Psi, \mathbf{h}_\alpha)$, of the first-order sensitivity $S_2(T^o, \Psi^0; \boldsymbol{\alpha}^0)$, in the form

$$DS_2(T^o, \Psi^0; \boldsymbol{\alpha}^0; h_T, h_\Psi, \mathbf{h}_\alpha) = \int_{-\ell/2}^{\ell/2} h_\Psi(z) \delta\left(z - \frac{\ell}{2}\right) dz, \tag{68}$$

As before, the function $h_\Psi$ in Eq. (68) is the solution of Eqs. (40) – (42). Following the same procedure as in the previous sub-section leads to the following expression for the partial-differential $DS_2(T^o, \Psi^0; \boldsymbol{\alpha}^0; h_T, h_\Psi, \mathbf{h}_\alpha)$:

$$\begin{aligned}DS_2(T^0; \boldsymbol{\alpha}^0; h_T, h_\Psi, \mathbf{h}_\alpha) &= \int_{-\ell/2}^{\ell/2} \left[\Psi_{21}^{(2)}(z) Q_1(T^0, \boldsymbol{\alpha}^0; \mathbf{h}_\alpha) + \Psi_{22}^{(2)}(z) Q_{12}(T^0, \boldsymbol{\alpha}^0; \mathbf{h}_\alpha)\right] dz \\ &\quad - q_1(T^0, \boldsymbol{\alpha}; \mathbf{h}_\alpha)\left[\Psi_{21}^{(2)}(z)\right]_{z=\ell/2} - (\delta T_a)\left[k(T^0)\frac{d\Psi_{21}^{(2)}(z)}{dz}\right]_{z=-\ell/2},\end{aligned} \tag{69}$$

where the 2$^{nd}$-level adjoint function $\boldsymbol{\Psi}_2^{(2)} \equiv \left(\Psi_{21}^{(2)}, \Psi_{22}^{(2)}\right)$ is the solution of the following *second-level adjoint sensitivity system (2$^{nd}$-LASS)*:

$$\begin{pmatrix} k(T^0)\dfrac{d^2[\ ]}{dz^2} & \dfrac{c^0 \delta(z-z_r)}{1+c^0 T^0(z)} \\ 0 & \dfrac{d^2}{dz^2}\{k(T^0)[\ ]\} \end{pmatrix} \begin{pmatrix} \Psi_{21}^{(2)}(z) \\ \Psi_{22}^{(2)}(z) \end{pmatrix} = \begin{pmatrix} 0 \\ \delta\left(z-\dfrac{\ell}{2}\right) \end{pmatrix}, \tag{70}$$

$$\left.\frac{d\Psi_{21}^{(2)}}{dz}\right|_{z=\ell/2} = 0, \quad at \ z = \ell/2, \tag{71}$$

$$\Psi_{21}^{(2)}(z) = 0, \quad at \ z = -\ell/2. \tag{72}$$

$$\left. \frac{d\left[ k\left(T^{0}\right)\Psi_{22}^{(2)} \right]}{dz} \right|_{z=\ell/2} = 0, \quad at \ z = \ell/2, \tag{73}$$

$$\Psi_{22}^{(2)}(z) = 0, \quad at \ z = -\ell/2. \tag{74}$$

Solving Eqs. (70) - (74) yields the following expressions for the components of the 2$^{nd}$-level adjoint function $\mathbf{\Psi}_{2}^{(2)} \equiv \left( \Psi_{21}^{(2)}, \Psi_{22}^{(2)} \right)$:

$$\Psi_{21}^{(2)}(z) = -C\left(T^{0}, z_{r}; \boldsymbol{\alpha}^{0}\right)\left( z_{r} + \frac{\ell}{2} \right)\left[ z + \frac{\ell}{2} - (z - z_{r})H(z - z_{r}) \right], \tag{75}$$

and

$$\Psi_{22}^{(2)}(z) = \frac{-1}{k\left[T^{0}(z)\right]}\left[ z + \frac{\ell}{2} - \left( z - \frac{\ell}{2} \right)H\left( z - \frac{\ell}{2} \right) \right]. \tag{76}$$

The second-order sensitivities $S_{2i} \triangleq \partial^{2}T(z_{r})/(\partial q \partial \alpha_{i})$ will have formally the same expressions as those shown Eqs. (58) – (62), except that the function $\mathbf{\Psi}_{1}^{(2)} \equiv \left( \Psi_{11}^{(2)}, \Psi_{12}^{(2)} \right)$ will be replaced by the function $\mathbf{\Psi}_{2}^{(2)} \equiv \left( \Psi_{21}^{(2)}, \Psi_{22}^{(2)} \right)$, to obtain:

$$\frac{\partial^{2}T(z_{r})}{\partial q \partial Q} = -\int_{-\ell/2}^{\ell/2} \Psi_{21}^{(2)}(z)\,dz, \tag{77}$$

$$\frac{\partial^{2}T(z_{r})}{\partial q^{2}} = \left[ \Psi_{21}^{(2)}(z) \right]_{z=\ell/2}, \tag{78}$$

$$\frac{\partial^{2}T(z_{r})}{\partial q \partial T_{a}} = -\left[ k\left(T^{0}\right) \frac{d\Psi_{21}^{(2)}(z)}{dz} \right]_{z=-\ell/2}, \tag{79}$$

$$\frac{\partial^{2}T(z_{r})}{\partial q \partial k_{0}} = \frac{1}{k_{0}^{0}}\left[ Q^{0}\int_{-\ell/2}^{\ell/2} \Psi_{21}^{(2)}(z)\,dz - \int_{-\ell/2}^{\ell/2} \Psi_{22}^{(2)}(z)\delta(z - z_{r})\,dz - q^{0}\left[ \Psi_{21}^{(2)}(z) \right]_{z=\ell/2} \right], \tag{80}$$

$$\frac{\partial^2 T(z_r)}{\partial q \partial c} = -k_0^0 \left\{ \frac{T_a^2}{2} \left[ \frac{d\Psi_{21}^{(2)}}{dz} \right]_{z=-\ell/2} + \frac{1}{2} \int_{-\ell/2}^{\ell/2} \left[ T^0(z) \right]^2 \frac{d^2 \Psi_{21}^{(2)}}{dz^2} dz + \frac{\Psi_{22}^{(2)}(z_r) T^0(z_r)}{k \left[ T^0(z_r) \right]} \right\}. \quad (81)$$

Replacing Eqs. (75) and (76) in the above expressions and carrying out the integrations over $z$ yields the following explicit expressions for the above 2$^{\text{nd}}$-order sensitivities:

$$\frac{\partial^2 T(z_r)}{\partial q \partial Q} = -\int_{-\ell/2}^{\ell/2} \Psi_{21}^{(2)}(z) dz = \frac{\partial^2 T(z_r)}{\partial Q \partial q} = \left[ \Psi_{11}^{(2)}(z) \right]_{z=\ell/2}$$
$$= \frac{1}{2} C(T^0, z_r; \boldsymbol{\alpha}^0) \left( \frac{3\ell^2}{4} + z_r \ell - z_r^2 \right) \left( z_r + \frac{\ell}{2} \right), \quad (82)$$

$$\frac{\partial^2 T(z_r)}{\partial q^2} = -C(T^0, z_r; \boldsymbol{\alpha}^0) \left( z_r + \frac{\ell}{2} \right)^2, \quad (83)$$

$$\frac{\partial^2 T(z_r)}{\partial q \partial T_a} = C(T^0, z_r; \boldsymbol{\alpha}^0) k(T_a^0) \left( z_r + \frac{\ell}{2} \right), \quad (84)$$

$$\frac{\partial^2 T(z_r)}{\partial q \partial k_0} = C(T^0, z_r; \boldsymbol{\alpha}^0) \left\{ \frac{1}{c^0} \left[ 1 + c^0 T^0(z_r) \right]^2 - \tau^0(z_r) \right\} \left( z_r + \frac{\ell}{2} \right), \quad (85)$$

$$\frac{\partial^2 T(z_r)}{\partial q \partial c} = C(T^0, z_r; \boldsymbol{\alpha}^0) \frac{k_0^0}{c^0} \left\{ c^0 \left[ T^0(z_r) \right]^2 + 2T^0(z_r) - \tau^0(z_r) - T_a^0 \right\} \left( z_r + \frac{\ell}{2} \right). \quad (86)$$

The symmetry of the second-order sensitivity $\partial^2 T(z_r)/\partial q \partial Q$ implies the equality between Eqs. (82) and (59) which, in turn, provides a stringent independent verification of the accuracy of computing the second-level adjoint functions $\boldsymbol{\Psi}_1^{(2)} \equiv \left( \Psi_{11}^{(2)}, \Psi_{12}^{(2)} \right)$ and $\boldsymbol{\Psi}_2^{(2)} \equiv \left( \Psi_{21}^{(2)}, \Psi_{22}^{(2)} \right)$. The 2$^{\text{nd}}$-order relative sensitivities $(S_{2i})_{rel} \triangleq \frac{\partial^2 T(z_r)}{\partial q \partial \alpha_i} \frac{q \alpha_i}{T(z_r)}$ are depicted in Figures 12 - 15 as functions of the arbitrary location $z_r \in [-\ell/2, \ell/2]$.

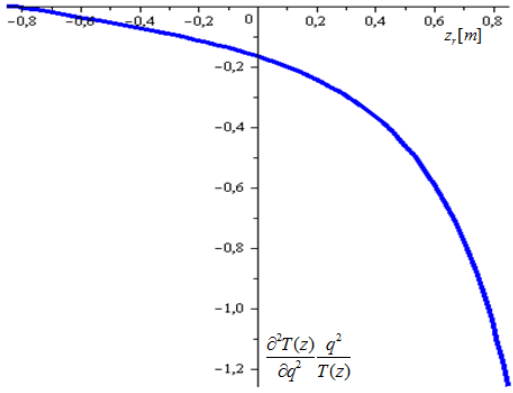
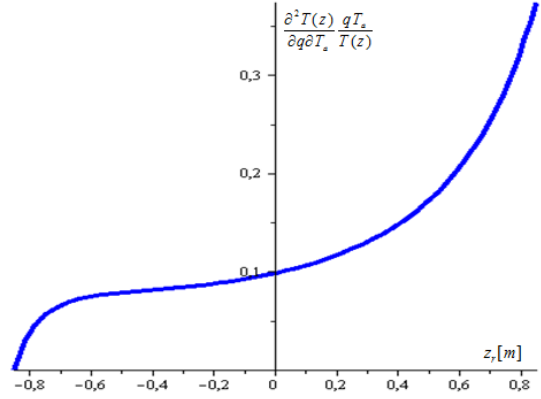

Fig. 12: Relative sensitivity $\dfrac{\partial^2 T(z_r)}{\partial q^2}\dfrac{q^2}{T(z_r)}$    Fig. 13: Relative sensitivity $\dfrac{\partial^2 T(z_r)}{\partial q \partial T_a}\dfrac{qT_a}{T(z_r)}$

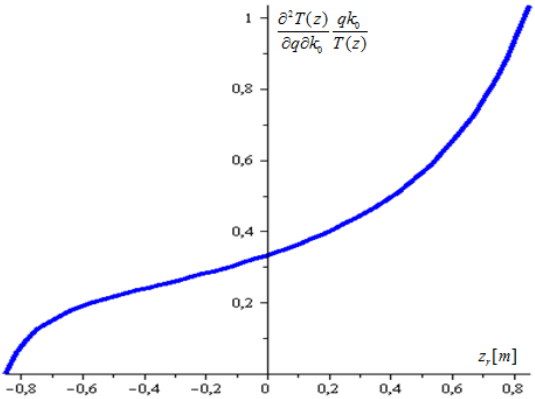
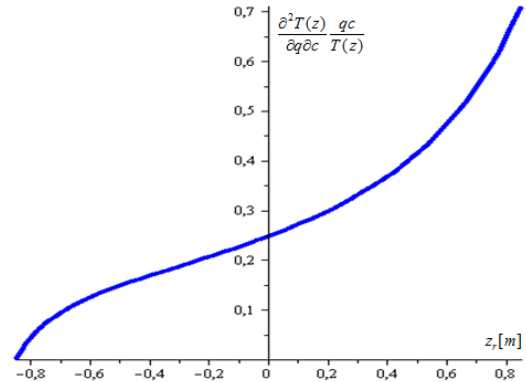

Fig. 14: Relative sensitivity $\dfrac{\partial^2 T(z_r)}{\partial q \partial k_0}\dfrac{qk_0}{T(z_r)}$    Fig. 15: Relative sensitivity $\dfrac{\partial^2 T(z_r)}{\partial q \partial c}\dfrac{qc}{T(z_r)}$

Figures 12 through 15 indicate, the the relative importance (magnitude) of the 2$^{nd}$-order relative sensitivities $(S_{2i})_{rel} \triangleq \dfrac{\partial^2 T(z_r)}{\partial q \partial \alpha_i}\dfrac{q\alpha_i}{T(z_r)}$ depends on the position $z_r$. Ranking them by the largest attainable maximum absolute values of these relative sensitivities, the importance of the model parameters is as follows: the heat source $Q$; the boundary heat flux $q$; the linear heat conductivity coefficient $k_0$; the nonlinear heat conductivity coefficient $c$; and the ambient temperature $T_a$.

### 3.2.3. Computation of the Second-Order Response Sensitivities $S_{3i} \triangleq \partial^2 T(z_r)/(\partial T_a \partial \alpha_i)$

Applying the definition shown in Eq. (37) to Eq. (29) yields the G-differential, $DS_3\left(T^o, \Psi^0; \boldsymbol{\alpha}^0; h_T, h_\Psi, \mathbf{h}_\alpha\right)$, of the first-order sensitivity $S_3\left(T^o, \Psi^0; \boldsymbol{\alpha}^0\right)$, in the form

$$
\begin{aligned}
& DS_3\left(T^o, \Psi^0; \boldsymbol{\alpha}^0; h_T, h_\Psi, \mathbf{h}_\alpha\right) = \\
& = \frac{d}{d\varepsilon}\left\{-\int_{-\ell/2}^{\ell/2} \left(k_0^0 + \varepsilon \delta k_0\right)\left[1 + \left(c^0 + \varepsilon \delta c\right)\left(T^0 + \varepsilon h_T\right)\right] \frac{d\left(\Psi^0 + \varepsilon h_\Psi\right)}{dz} \delta\left(z + \frac{\ell}{2}\right) dz\right\}_{\varepsilon=0} \\
& = \left[DS_3\left(T^o, \Psi^0; \boldsymbol{\alpha}^0; h_T, h_\Psi, \mathbf{h}_\alpha\right)\right]_{direct} + \left[DS_3\left(T^o, \Psi^0; \boldsymbol{\alpha}^0; h_T, h_\Psi, \mathbf{h}_\alpha\right)\right]_{indirect},
\end{aligned}
\quad (87)
$$

where the "direct-effect term" is defined as

$$
\begin{aligned}
& \left[DS_3\left(T^o, \Psi^0; \boldsymbol{\alpha}^0; h_T, h_\Psi, \mathbf{h}_\alpha\right)\right]_{direct} \equiv \\
& = -\int_{-\ell/2}^{\ell/2}\left\{\left(\delta k_0\right)\left[1 + c^0 T^0(z)\right] + \left(\delta c\right) k_0^0 T^0(z) + c^0 k_0^0 h_T(z)\right\}\frac{d\Psi^0(z)}{dz}\delta\left(z + \frac{\ell}{2}\right)dz,
\end{aligned}
\quad (88)
$$

while the "indirect-effect term" is defined as

$$
\left[DS_3\left(T^o, \Psi^0; \boldsymbol{\alpha}^0; h_T, h_\Psi, \mathbf{h}_\alpha\right)\right]_{indirect} \equiv -\int_{-\ell/2}^{\ell/2} k\left[T^0(z)\right]\frac{dh_\Psi(z)}{dz}\delta\left(z + \frac{\ell}{2}\right)dz. \quad (89)
$$

The "direct-effect term" defined in Eq. (88) can be evaluated immediately by noting from Eq. (23) that

$$
\left[\frac{d\Psi(z)}{dz}\right]_{z=-\ell/2} = -\frac{1}{k\left[T^0(z_r)\right]}, \quad (90)
$$

and by using the above result in Eq. (88) to obtain

$$
\left[DS_3\left(T^o, \Psi^0; \boldsymbol{\alpha}^0; h_T, h_\Psi, \mathbf{h}_\alpha\right)\right]_{direct} = \frac{1}{k\left[T^0(z_r)\right]}\left[\left(\delta k_0\right)\left(1 + c^0 T_a^0\right) + \left(\delta c\right)k_0^0 T_a^0 + \left(\delta T_a\right)c^0 k_0^0\right]. \quad (91)
$$

On the other hand, the "indirect-effect term" defined in Eq. (89) needs to be evaluated by constructing the corresponding *2$^{nd}$-LASS* for a 2$^{nd}$-level adjoint function $\mathbf{\Psi}_3^{(2)} \equiv \left( \Psi_{31}^{(2)}, \Psi_{32}^{(2)} \right)$ by following the general general principles of the *2$^{nd}$-ASAM* presented in PART I [1]. Applying these principles leads to the following expression for the "indirect-effect term" defined in Eq. (89):

$$\left[ DS_3 \left( T^o, \Psi^0; \boldsymbol{\alpha}^0; h_T, h_\Psi, \mathbf{h}_\alpha \right) \right]_{indirect}$$
$$= \int_{-\ell/2}^{\ell/2} \left[ \Psi_{31}^{(2)}(z) Q_1 \left( T^0, \boldsymbol{\alpha}^0; \mathbf{h}_\alpha \right) + \Psi_{32}^{(2)}(z) Q_{12} \left( T^0, \boldsymbol{\alpha}^0; \mathbf{h}_\alpha \right) \right] dz \qquad (92)$$
$$- q_1 \left( T^0, \boldsymbol{\alpha}; \mathbf{h}_\alpha \right) \left[ \Psi_{31}^{(2)}(z) \right]_{z=\ell/2} - \left( \delta T_a \right) \left[ k \left( T^0 \right) \frac{d \Psi_{31}^{(2)}(z)}{dz} \right]_{z=-\ell/2},$$

where the 2$^{nd}$-level adjoint function $\mathbf{\Psi}_3^{(2)} \equiv \left( \Psi_{31}^{(2)}, \Psi_{32}^{(2)} \right)$ is the solution of the following *2$^{nd}$-LASS*:

$$\begin{pmatrix} k(T^0) \dfrac{d^2 [\ ]}{dz^2} & \dfrac{c^0 \delta(z - z_r)}{1 + c^0 T^0(z)} \\ 0 & \dfrac{d^2}{dz^2} \left\{ k(T^0) [\ ] \right\} \end{pmatrix} \begin{pmatrix} \Psi_{31}^{(2)}(z) \\ \Psi_{32}^{(2)}(z) \end{pmatrix} = \begin{pmatrix} 0 \\ 0 \end{pmatrix}, \qquad (93)$$

$$\left. \frac{d \Psi_{31}^{(2)}}{dz} \right|_{z=\ell/2} = 0, \quad at \ z = \ell/2, \qquad (94)$$

$$\Psi_{31}^{(2)}(z) = 0, \quad at \ z = -\ell/2. \qquad (95)$$

$$\left. \frac{d \left[ k(T^0) \Psi_{32}^{(2)} \right]}{dz} \right|_{z=\ell/2} = 0, \quad at \ z = \ell/2, \qquad (96)$$

$$\Psi_{32}^{(2)}(z) = 1, \quad at \ z = -\ell/2. \qquad (97)$$

Solving Eqs. (93) - (97) yields the following expressions for the components of the 2nd-level adjoint function $\mathbf{\Psi}_3^{(2)} \equiv \left( \Psi_{31}^{(2)}, \Psi_{32}^{(2)} \right)$:

$$\Psi_{31}^{(2)}(z) = C\left(T^0, z_r; \boldsymbol{\alpha}^0\right) k(T_a) \left[ z + \frac{\ell}{2} - (z - z_r) H(z - z_r) \right], \tag{98}$$

and

$$\Psi_{32}^{(2)}(z) = \frac{k(T_a)}{k\left[T^0(z)\right]}. \tag{99}$$

Adding Eqs. (88) and (92) and identifying the coefficients multiplying the respective parameter variations yields the following expressions for the 2nd-order sensitivities $S_{3i} \triangleq \partial^2 T(z_r) / (\partial T_a \partial \alpha_i)$:

$$\begin{aligned}
\frac{\partial^2 T(z_r)}{\partial T_a \partial Q} &= -\int_{-\ell/2}^{\ell/2} \Psi_{31}^{(2)}(z) dz = \frac{\partial^2 T(z_r)}{\partial Q \partial T_a} = -\left[ k(T^0) \frac{d\Psi_{11}^{(2)}(z)}{dz} \right]_{z=-\ell/2} \\
&= -C\left(T^0, z_r; \boldsymbol{\alpha}^0\right) \frac{k(T_a^0)}{2} \left( \frac{3\ell^2}{4} + z_r \ell - z_r^2 \right),
\end{aligned} \tag{100}$$

$$\begin{aligned}
\frac{\partial^2 T(z_r)}{\partial T_a \partial q} &= \left[ \Psi_{31}^{(2)}(z) \right]_{z=\ell/2} = \frac{\partial^2 T(z_r)}{\partial q \partial T_a} = -\left[ k(T^0) \frac{d\Psi_{21}^{(2)}(z)}{dz} \right]_{z=-\ell/2} \\
&= C\left(T^0, z_r; \boldsymbol{\alpha}^0\right) k(T_a^0) \left( z_r + \frac{\ell}{2} \right),
\end{aligned} \tag{101}$$

$$\begin{aligned}
\frac{\partial^2 T(z_r)}{\partial T_a^2} &= -\left[ c^0 k_0^0 \frac{d\Psi^0(z)}{dz} + k(T^0) \frac{d\Psi_{31}^{(2)}(z)}{dz} \right]_{z=-\ell/2} \\
&= C\left(T^0, z_r; \boldsymbol{\alpha}^0\right) \left\{ k^2\left[T^0(z_r)\right] - k^2\left(T_a^0\right) \right\},
\end{aligned} \tag{102}$$

$$\begin{aligned}
\frac{\partial^2 T(z_r)}{\partial T_a \partial k_0} &= \frac{1 + c^0 T_a^0}{k\left[T^0(z_r)\right]} + \frac{Q^0}{k_0^0} \int_{-\ell/2}^{\ell/2} \Psi_{31}^{(2)}(z) dz \\
&\quad - \frac{1}{k_0^0} \int_{-\ell/2}^{\ell/2} \Psi_{32}^{(2)}(z) \delta(z - z_r) dz - \frac{q^0}{k_0^0} \left[ \Psi_{31}^{(2)}(z) \right]_{z=\ell/2} \\
&= C\left(T^0, z_r; \boldsymbol{\alpha}^0\right) k\left(T_a^0\right) \tau^0(z_r),
\end{aligned} \tag{103}$$

$$\frac{\partial^2 T(z_r)}{\partial T_a \partial c} = \frac{k_0^0 T_a^0}{k\left[T^0(z_r)\right]} + C\left(T^0, z_r; \boldsymbol{\alpha}^0\right) \frac{k_0^0 k\left(T_a^0\right)}{c^0} \left\{ \tau^0(z_r) + T_a^0 - c^0 \left[T^0(z_r)\right]^2 - 2T^0(z_r) \right\}. \quad (104)$$

The symmetry of the second-order sensitivity $\partial^2 T(z_r)/\partial T_a \partial Q$ implies the equality between the relations expressed in Eq. (100) and (65). This equality provides an independent verification of the correctness of the respective expressions as well as a verification of the solution accuracy of computing the second-level adjoint functions $\boldsymbol{\Psi}_1^{(2)} \equiv \left(\Psi_{11}^{(2)}, \Psi_{12}^{(2)}\right)$ and $\boldsymbol{\Psi}_3^{(2)} \equiv \left(\Psi_{31}^{(2)}, \Psi_{32}^{(2)}\right)$. Furthermore, the symmetry of the second-order sensitivity $\partial^2 T(z_r)/\partial T_a \partial q$ implies the equality between the relations expressed in Eq. (101) and (79). This equality provides an independent verification of the correctness of the respective expressions as well as a verification of the solution accuracy of computing the second-level adjoint functions $\boldsymbol{\Psi}_2^{(2)} \equiv \left(\Psi_{21}^{(2)}, \Psi_{22}^{(2)}\right)$ and $\boldsymbol{\Psi}_3^{(2)} \equiv \left(\Psi_{31}^{(2)}, \Psi_{32}^{(2)}\right)$. The 2$^{nd}$-order relative sensitivities $(S_{3i})_{rel} \triangleq \frac{\partial^2 T(z_r)}{\partial T_a \partial \alpha_i} \frac{T_a \alpha_i}{T(z_r)}$ are depicted in Figures 16 - 18 as functions of the arbitrary location $z_r \in [-\ell/2, \ell/2]$.

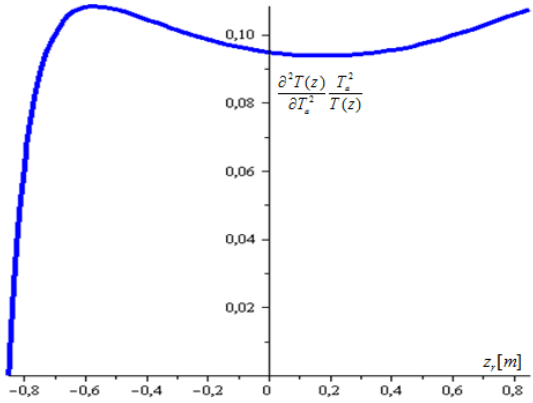 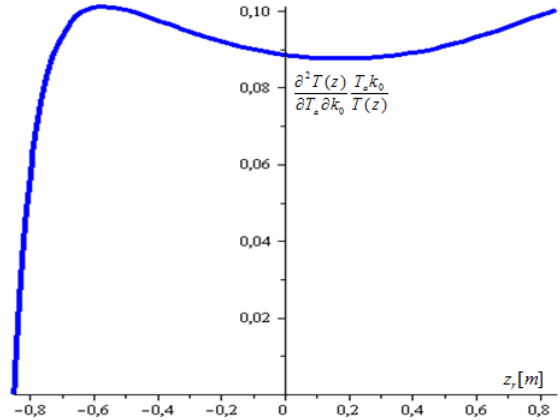

Fig. 16: Relative sensitivity $\dfrac{\partial^2 T(z_r)}{\partial T_a^2} \dfrac{T_a^2}{T(z_r)}$      Fig. 17: Relative sensitivity $\dfrac{\partial^2 T(z_r)}{\partial T_a \partial k_0} \dfrac{T_a k_0}{T(z_r)}$

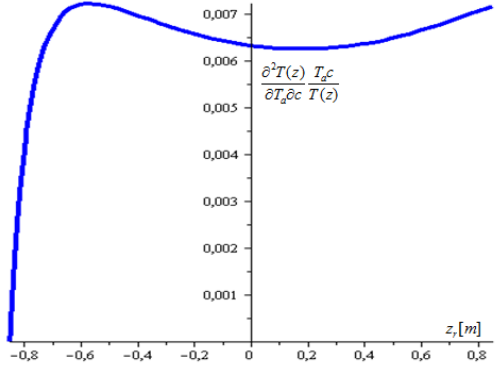

Fig. 18: Relative sensitivity $\dfrac{\partial^2 T(z_r)}{\partial T_a \partial c} \dfrac{T_a c}{T(z_r)}$

As Figures 16 through 18 indicate, the the relative importance (magnitude) of the 2$^{nd}$-order relative sensitivities $(S_{3i})_{rel} \triangleq \dfrac{\partial^2 T(z_r)}{\partial T_a \partial \alpha_i} \dfrac{T_a \alpha_i}{T(z_r)}$ depends on the position $z_r$. Ranking them by the largest attainable maximum absolute values of these relative sensitivities, the importance of the model parameters is as follows: the boundary heat flux $q$, as previously depicted in Figure 13; and the heat source $Q$, as previously depicted in Figure 9. On the other hand, Figures 16 -18 indicate that the remaining parameters (namely the linear heat conductivity coefficient $k_0$, the ambient temperature $T_a$, and the nonlinear heat conductivity coefficient $c$) are much less important.

*3.2.4. Computation of the Second-Order Response Sensitivities $S_{4i} \triangleq \partial^2 T(z_r)/(\partial k_0 \partial \alpha_i)$*

Applying the definition shown in Eq. (37) to Eq. (30) yields the G-differential, $DS_4(T^o, \Psi^0; \boldsymbol{\alpha}^0; h_T, h_\Psi, \mathbf{h}_\alpha)$, of the first-order sensitivity $S_4(T^o, \Psi^0; \boldsymbol{\alpha}^0)$, in the form

$$DS_4(T^o, \Psi^0; \boldsymbol{\alpha}^0; h_T, h_\Psi, \mathbf{h}_\alpha) = \left[DS_4(T^o, \Psi^0; \boldsymbol{\alpha}^0; h_T, h_\Psi, \mathbf{h}_\alpha)\right]_{direct} \\ + \left[DS_4(T^o, \Psi^0; \boldsymbol{\alpha}^0; h_T, h_\Psi, \mathbf{h}_\alpha)\right]_{indirect}, \qquad (105)$$

where the "direct effect" term is defined as $\left[DS_4(T^o, \Psi^0; \boldsymbol{\alpha}^0; h_T, h_\Psi, \mathbf{h}_\alpha)\right]_{direct}$

$$\left[ DS_4 \left(T^o, \Psi^0; \boldsymbol{\alpha}^0; h_T, h_\Psi, \mathbf{h}_\alpha \right) \right]_{direct} \equiv \left[ \frac{(\delta Q)}{k_0^0} - (\delta k_0) \frac{Q^0}{\left(k_0^0\right)^2} \right] \int_{-\ell/2}^{\ell/2} \Psi(z) dz$$
$$- \left[ \frac{(\delta q)}{k_0^0} - (\delta k_0) \frac{q^0}{\left(k_0^0\right)^2} \right] \Psi(\ell/2), \quad (106)$$

while the "indirect effect" $\left[ DS_4 \left(T^o, \Psi^0; \boldsymbol{\alpha}^0; h_T, h_\Psi, \mathbf{h}_\alpha \right) \right]_{indirect}$ term is defined as

$$\left[ DS_4 \left(T^o, \Psi^0; \boldsymbol{\alpha}^0; h_T, h_\Psi, \mathbf{h}_\alpha \right) \right]_{indirect} \equiv \frac{Q^0}{k_0^0} \int_{-\ell/2}^{\ell/2} h_\Psi(z) dz - \frac{q^0}{k_0^0} h_\Psi(\ell/2). \quad (107)$$

The "direct effect" term defined in Eq. (106) can be evaluated by using Eq. (23) to obtain:

$$\left[ DS_4 \left(T^o, \Psi^0; \boldsymbol{\alpha}^0; h_T, h_\Psi, \mathbf{h}_\alpha \right) \right]_{direct} = \left[ (\delta k_0) \frac{Q^0}{\left(k_0^0\right)^2} - \frac{(\delta Q)}{k_0^0} \right] \frac{1}{2k\left[T^0(z_r)\right]} \left( \frac{3\ell^2}{4} + z_r \ell - z_r^2 \right)$$
$$+ \left[ \frac{(\delta q)}{k_0^0} - (\delta k_0) \frac{q^0}{\left(k_0^0\right)^2} \right] \frac{z_r + \ell/2}{k\left[T^0(z_r)\right]}. \quad (108)$$

The "indirect-effect term" defined in Eq. (107) needs to be evaluated by constructing the corresponding *2nd-LASS* for a 2nd-level adjoint function $\boldsymbol{\Psi}_4^{(2)} \equiv \left(\Psi_{41}^{(2)}, 0\right)$ by following the general general principles of the *2nd-ASAM* presented in PART I [1]. Applying these principles leads to the following expression for the "indirect-effect term" defined in Eq. (107):

$$\left[ DS_4 \left(T^o, \Psi^0; \boldsymbol{\alpha}^0; h_T, h_\Psi, \mathbf{h}_\alpha \right) \right]_{indirect} = \int_{-\ell/2}^{\ell/2} \left[ \Psi_{41}^{(2)}(z) Q_1 \left(T^0, \boldsymbol{\alpha}^0; \mathbf{h}_\alpha \right) + \Psi_{42}^{(2)}(z) Q_{12} \left(T^0, \boldsymbol{\alpha}^0; \mathbf{h}_\alpha \right) \right] dz$$
$$- q_1 \left(T^0, \boldsymbol{\alpha}; \mathbf{h}_\alpha \right) \left[ \Psi_{41}^{(2)}(z) \right]_{z=\ell/2} - (\delta T_a) \left[ k(T^0) \frac{d\Psi_{41}^{(2)}(z)}{dz} \right]_{z=-\ell/2},$$
$$(109)$$

where the 2nd-level adjoint function $\boldsymbol{\Psi}_4^{(2)} \equiv \left(\Psi_{41}^{(2)}, \Psi_{42}^{(2)}\right)$ is the solution of the following *2nd-LASS*:

$$\begin{pmatrix} k(T^0)\dfrac{d^2[\ ]}{dz^2} & \dfrac{c^0\delta(z-z_r)}{1+c^0T^0(z)} \\ 0 & \dfrac{d^2}{dz^2}\{k(T^0)[\ ]\} \end{pmatrix} \begin{pmatrix} \Psi_{41}^{(2)}(z) \\ \Psi_{42}^{(2)}(z) \end{pmatrix} = \begin{pmatrix} 0 \\ \dfrac{Q^0}{k_0^0} - \dfrac{q^0}{k_0^0}\delta\left(z-\dfrac{\ell}{2}\right) \end{pmatrix}, \quad (110)$$

$$\left.\frac{d\Psi_{41}^{(2)}}{dz}\right|_{z=\ell/2} = 0, \quad at\ z=\ell/2, \tag{111}$$

$$\Psi_{41}^{(2)}(z) = 0, \quad at\ z=-\ell/2. \tag{112}$$

$$\left.\frac{d\left[k(T^0)\Psi_{42}^{(2)}\right]}{dz}\right|_{z=\ell/2} = 0, \quad at\ z=\ell/2, \tag{113}$$

$$\Psi_{42}^{(2)}(z) = 0, \quad at\ z=-\ell/2. \tag{114}$$

Solving Eqs (110) – (114) yields the following expressions for the 2$^{nd}$-level adjoint function $\boldsymbol{\Psi}_4^{(2)} \equiv \left(\Psi_{41}^{(2)}, \Psi_{42}^{(2)}\right)$:

$$\Psi_{41}^{(2)}(z) = -C\left(T^0, z_r; \boldsymbol{\alpha}^0\right)\tau(z_r)\left(\frac{3\ell^2}{4} + z_r\ell - z_r^2\right)\left[z + \frac{\ell}{2} - (z-z_r)H(z-z_r)\right], \tag{115}$$

and

$$\Psi_{42}^{(2)}(z) = \frac{1}{k\left[T^0(z)\right]}\left\{-\frac{Q^0}{2k_0^0}\left(\frac{3\ell^2}{4} + z_r\ell - z_r^2\right) + \frac{q^0}{k_0^0}\left[z + \frac{\ell}{2} - \left(z-\frac{\ell}{2}\right)H\left(z-\frac{\ell}{2}\right)\right]\right\}, \tag{116}$$

Replacing Eq. (115) and (116) in Eq. (109), carrying out the integrations over $z$ and the ensuing algebra yields the following expression for the "indirect-effect term" $\left[DS_5\left(T^o, \Psi^0; \boldsymbol{\alpha}^0; h_T, h_\Psi, \mathbf{h}_\alpha\right)\right]_{indirect}$:

$$\left[DS_4\left(T^o,\Psi^0;\boldsymbol{\alpha}^0;h_T,h_\Psi,\mathbf{h}_\alpha\right)\right]_{indirect} = \left[\left(\delta Q\right)-\left(\delta k_0\right)\frac{Q^0}{k_0^0}\right]\frac{c^0 k_0^0 \tau(z_r)}{2k^3\left[T^0(z_r)\right]}\left(\frac{3\ell^2}{4}+z_r\ell-z_r^2\right)$$

$$+\left(\delta c\right)\frac{c^0\left(k_0^0\right)^2 \tau(z_r)}{2k^3\left[T^0(z_r)\right]}\left\{T_a^2-\left[T^0(z_r)\right]^2\right\}+\left[\frac{(\delta k_0)}{k_0^0}+(\delta c)\frac{T^0(z_r)}{1+c^0 T^0(z_r)}\right]\frac{\tau(z_r)}{k\left[T^0(z_r)\right]} \quad (117)$$

$$+(\delta T_a) k(T_a)\frac{c^0 k_0^0 \tau(z_r)}{k^3\left[T^0(z_r)\right]}+\left[(\delta k_0)\frac{q^0}{k_0^0}-(\delta q)\right]\frac{c^0 k_0^0 \tau(z_r)}{k^3\left[T^0(z_r)\right]}\left(z_r+\frac{\ell}{2}\right).$$

Summing the expression for the "indirect effect term" given in Eq. (117) with the expression of the "direct effect term" given in Eq. (108) and identifying the coefficients of the respective parameter variations yields the following expression for the 2$^{nd}$-order derivatives $S_{4i} \triangleq \partial^2 T(z_r)/(\partial k_0 \partial \alpha_i)$:

$$\frac{\partial^2 T(z_r)}{\partial k_0 \partial Q} = C(T^0,z_r;\boldsymbol{\alpha}^0)\frac{1}{2}\left\{\tau^0(z_r)-\frac{1}{c^0}\left[1+c^0 T^0(z_r)\right]^2\right\}\left(\frac{3\ell^2}{4}+z_r\ell-z_r^2\right), \quad (118)$$

$$\frac{\partial^2 T(z_r)}{\partial k_0 \partial q} = C(T^0,z_r;\boldsymbol{\alpha}^0)\left\{\frac{1}{c^0}\left[1+c^0 T^0(z_r)\right]^2-\tau^0(z_r)\right\}\left(z_r+\frac{\ell}{2}\right), \quad (119)$$

$$\frac{\partial^2 T(z_r)}{\partial k_0 \partial T_a} = C(T^0,z_r;\boldsymbol{\alpha}^0) k(T_a^0) \tau^0(z_r), \quad (120)$$

$$\frac{\partial^2 T(z_r)}{\partial k_0^2} = \tau^0(z_r)\left\{\frac{2}{k_0^0 k\left[T^0(z_r)\right]}-C(T^0,z_r;\boldsymbol{\alpha}^0)\tau^0(z_r)\right\}, \quad (121)$$

$$\frac{\partial^2 T(z_r)}{\partial k_0 \partial c} = C(T^0,z_r;\boldsymbol{\alpha}^0)\frac{k_0^0 \tau^0(z_r)}{c^0}\left\{c^0\left[T^0(z_r)\right]^2+2T^0(z_r)-\tau^0(z_r)-T_a^0\right\}, \quad (122)$$

The symmetry of the second-order sensitivity $\partial^2 T(z_r)/\partial k_0 \partial Q$ implies the equality between the relations expressed in Eq. (118) and (66). This equality provides an independent verification of the correctness of the respective expressions as well as a verification of the solution accuracy of computing the second-level adjoint functions $\boldsymbol{\Psi}_1^{(2)} \equiv \left(\Psi_{11}^{(2)}, \Psi_{12}^{(2)}\right)$ and $\boldsymbol{\Psi}_4^{(2)} \equiv \left(\Psi_{41}^{(2)}, \Psi_{42}^{(2)}\right)$. Furthermore, the symmetry of the second-order sensitivity $\partial^2 T(z_r)/\partial k_0 \partial q$ implies the equality between the relations expressed in Eq. (119) and (85). This equality provides an independent verification of the correctness of the respective expressions as well as a verification of the solution accuracy of computing the second-level adjoint functions $\mathbf{\Psi}_2^{(2)} \equiv \left( \Psi_{21}^{(2)}, \Psi_{22}^{(2)} \right)$ and $\mathbf{\Psi}_4^{(2)} \equiv \left( \Psi_{41}^{(2)}, \Psi_{42}^{(2)} \right)$. Finally, the symmetry of the second-order sensitivity $\partial^2 T(z_r)/\partial k_0 \partial T_a$ implies the equality between the relations expressed in Eq. (120) and (103). This equality provides an independent verification of the correctness of the respective expressions as well as a verification of the solution accuracy of computing the second-level adjoint functions $\mathbf{\Psi}_3^{(2)} \equiv \left( \Psi_{31}^{(2)}, \Psi_{32}^{(2)} \right)$ and $\mathbf{\Psi}_4^{(2)} \equiv \left( \Psi_{41}^{(2)}, \Psi_{42}^{(2)} \right)$. The 2nd-order relative sensitivities $\dfrac{\partial^2 T(z_r)}{\partial k_0^2} \dfrac{k_0^2}{T(z_r)}$ and $\dfrac{\partial^2 T(z_r)}{\partial k_0 \partial c} \dfrac{k_0 c}{T(z_r)}$ are depicted in Figures 19 and 20 as functions of the arbitrary location $z_r \in [-\ell/2, \ell/2]$.

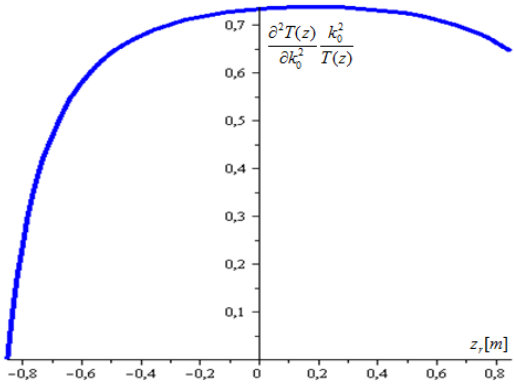
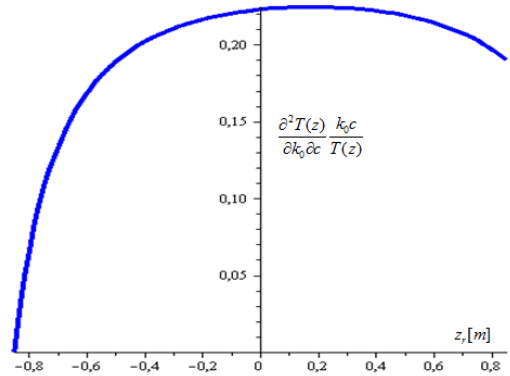

Fig. 19: Relative sensitivity $\dfrac{\partial^2 T(z_r)}{\partial k_0^2} \dfrac{k_0^2}{T(z_r)}$    Fig. 20: Relative sensitivity $\dfrac{\partial^2 T(z_r)}{\partial k_0 \partial c} \dfrac{k_0 c}{T(z_r)}$

*3.2.5. Computation of the Second-Order Response Sensitivities $S_{5i} \triangleq \partial^2 T(z_r)/(\partial c \partial \alpha_i)$*

Applying the definition shown in Eq. (37) to Eq. (31) yields the G-differential, $DS_5 \left( T^o, \Psi^0; \boldsymbol{\alpha}^0; h_T, h_\Psi, \mathbf{h}_\alpha \right)$, of the first-order sensitivity $S_5 \left( T^o, \Psi^0; \boldsymbol{\alpha}^0 \right)$, in the form

$$DS_5\left(T^o,\Psi^0;\boldsymbol{\alpha}^0;h_T,h_\Psi,\mathbf{h}_\alpha\right)=\frac{1}{2}\frac{d}{d\varepsilon}\left\{\int_{-\ell/2}^{\ell/2}\frac{\left(T_a^0+\varepsilon\delta T_a\right)^2-\left(T^0+\varepsilon h_T\right)^2}{\left[1+\left(c^0+\varepsilon\delta c\right)\left(T^0+\varepsilon h_T\right)\right]}\delta(z-z_r)dz\right\}_{\varepsilon=0} \quad (123)$$

$$=\left[DS_5\left(T^o,\Psi^0;\boldsymbol{\alpha}^0;h_T,h_\Psi,\mathbf{h}_\alpha\right)\right]_{direct}+\left[DS_5\left(T^o,\Psi^0;\boldsymbol{\alpha}^0;h_T,h_\Psi,\mathbf{h}_\alpha\right)\right]_{indirect},$$

where

$$\left[DS_5\left(T^o,\Psi^0;\boldsymbol{\alpha}^0;h_T,h_\Psi,\mathbf{h}_\alpha\right)\right]_{direct}\equiv\left(\delta T_a\right)\frac{T_a^0}{1+c^0T^0(z_r)}+\left(\delta c\right)\frac{\left[T^0(z_r)\right]^2-\left(T_a^0\right)^2}{2\left[1+c^0T^0(z_r)\right]^2}T^0(z_r), \quad (124)$$

and

$$\left[DS_5\left(T^o,\Psi^0;\boldsymbol{\alpha}^0;h_T,h_\Psi,\mathbf{h}_\alpha\right)\right]_{indirect}\equiv-\int_{-\ell/2}^{\ell/2}\frac{c^0\left(T_a^0\right)^2+2T^0(z)+c^0\left[T^0(z)\right]^2}{2\left[1+c^0T^0(z)\right]^2}h_T(z)\delta(z-z_r)dz. \quad (125)$$

The "indirect-effect term" defined in Eq. (125) needs to be evaluated by constructing the corresponding *2$^{nd}$-LASS* for a 2$^{nd}$-level adjoint function $\boldsymbol{\Psi}_5^{(2)}\equiv\left(\Psi_{51}^{(2)},0\right)$ by following the general general principles of the *2$^{nd}$-ASAM* presented in PART I [1]. Applying these principles leads to the following expression for the "indirect-effect term" defined in Eq. (125):

$$\left[DS_5\left(T^o,\Psi^0;\boldsymbol{\alpha}^0;h_T,h_\Psi,\mathbf{h}_\alpha\right)\right]_{indirect}=\int_{-\ell/2}^{\ell/2}\Psi_{51}^{(2)}(z)Q_1\left(T^0,\boldsymbol{\alpha}^0;\mathbf{h}_\alpha\right)dz$$
$$-q_1\left(T^0,\boldsymbol{\alpha};\mathbf{h}_\alpha\right)\left[\Psi_{51}^{(2)}(z)\right]_{z=\ell/2}-\left(\delta T_a\right)\left[k\left(T^0\right)\frac{d\Psi_{51}^{(2)}(z)}{dz}\right]_{z=-\ell/2}, \quad (126)$$

where the 2$^{nd}$-level adjoint function $\boldsymbol{\Psi}_5^{(2)}\equiv\left(\Psi_{51}^{(2)},0\right)$ is the solution of the following *2$^{nd}$-LASS*:

$$\frac{d^2\Psi_{51}^{(2)}(z)}{dz^2}=-\left(k_0^0\right)^2\frac{c^0\left(T_a^0\right)^2+2T^0(z)+c^0\left[T^0(z)\right]^2}{2k^3\left[T^0(z)\right]}\delta(z-z_r),\ z\in(-\ell/2,\ell/2), \quad (127)$$

$$\left.\frac{d\Psi_{51}^{(2)}}{dz}\right|_{z=\ell/2}=0,\ at\ z=\ell/2, \quad (128)$$

$$\Psi_{51}^{(2)}(z) = 0, \quad at \ z = -\ell/2. \tag{129}$$

Note that the above $2^{nd}$-LASS has the same form as the $1^{st}$-LASS [cf. Eqs. (20)-(22)], except for a different source term in Eq. (127). Therefore, the solution of Eqs. (127) – (29) can be written down by simply modifying appropriately the solution of the $1^{st}$-level adjoint function $\Psi(z)$ in Eq. (23) to obtain

$$\Psi_{51}^{(2)}(z) = -C(T^0, z_r; \boldsymbol{\alpha}^0) \frac{k_0^0}{2c^0} \left\{ c^0 (T_a^0)^2 + 2T^0(z_r) + c^0 \left[ T^0(z_r) \right]^2 \right\} \left[ (z - z_r) H(z - z_r) - z - \ell/2 \right]. \tag{130}$$

Replacing Eq. (130) in Eq. (126), carrying out the integrations over $z$ and the ensuing algebra yields the following expression for the "indirect-effect term" $\left[ DS_5(T^o, \Psi^0; \boldsymbol{\alpha}^0; h_T, h_\Psi, \mathbf{h}_\alpha) \right]_{indirect}$ :

$$\left[ DS_5(T^o, \Psi^0; \boldsymbol{\alpha}^0; h_T, h_\Psi, \mathbf{h}_\alpha) \right]_{indirect} =$$
$$-C(T^0, z_r; \boldsymbol{\alpha}^0) \frac{k_0^0}{2c^0} \left\{ c^0 (T_a^0)^2 + 2T^0(z_r) + c^0 \left[ T^0(z_r) \right]^2 \right\} \left\{ (\delta Q) \frac{1}{2} \left( \frac{3\ell^2}{4} + z_r \ell - z_r^2 \right) - (\delta q) \left( z_r + \frac{\ell}{2} \right) \right.$$
$$+ (\delta T_a) k (T_a) - (\delta k_0) \left[ \frac{Q^0}{2k_0^0} \left( \frac{3\ell^2}{4} + z_r \ell - z_r^2 \right) - \frac{q^0}{k_0^0} \left( z_r + \frac{\ell}{2} \right) \right] + (\delta c) \frac{(T_a^0)^2 - \left[ T^0(z_r) \right]^2}{2} \right\}. \tag{131}$$

Adding Eqs. (131) and (124), and identifying the coefficients multiplying the respective parameter variations yields the following expressions for the $2^{nd}$-order sensitivities $S_{5i} \triangleq \partial^2 T(z_r)/(\partial c \partial \alpha_i)$:

$$\frac{\partial^2 T(z_r)}{\partial c \partial Q} = C(T^0, z_r; \boldsymbol{\alpha}^0) \frac{k_0^0}{2c^0} \left\{ T_a^0 + \tau^0(z_r) - 2T^0(z_r) - c^0 \left[ T^0(z_r) \right]^2 \right\} \left( \frac{3\ell^2}{4} + z_r \ell - z_r^2 \right), \tag{132}$$

$$\frac{\partial^2 T(z_r)}{\partial c \partial q} = C(T^0, z_r; \boldsymbol{\alpha}^0) \frac{k_0^0}{c^0} \left\{ c^0 \left[ T^0(z_r) \right]^2 + 2T^0(z_r) - \tau^0(z_r) - T_a^0 \right\} \left( z_r + \frac{\ell}{2} \right), \tag{133}$$

$$\frac{\partial^2 T(z_r)}{\partial c \partial T_a} = \frac{k_0^0 T_a^0}{k[T^0(z_r)]} + C(T^0, z_r; \boldsymbol{\alpha}^0) \frac{k_0^0 k(T_a^0)}{c^0} \left\{ \tau^0(z_r) + T_a^0 - c^0 \left[T^0(z_r)\right]^2 - 2T^0(z_r) \right\}, \quad (134)$$

$$\frac{\partial^2 T(z_r)}{\partial c \partial k_0} = C(T^0, z_r; \boldsymbol{\alpha}^0) \frac{k_0^0 \tau^0(z_r)}{c^0} \left\{ c^0 \left[T^0(z_r)\right]^2 + 2T^0(z_r) - \tau^0(z_r) - T_a^0 \right\}, \quad (135)$$

$$\frac{\partial^2 T(z_r)}{\partial c^2} = C(T^0, z_r; \boldsymbol{\alpha}^0) \left(\frac{k_0^0}{c^0}\right)^2 \left[T_a^0 + \tau^0(z_r) - T^0(z_r)\right]$$
$$\times \left\{ c^0 \left[T^0(z_r)\right]^2 + c^0 \left(T_a^0\right)^2 + T^0(z_r) + \tau^0(z_r) + T_a^0 \right\}, \quad (136)$$

The symmetry of the second-order sensitivity $\partial^2 T(z_r)/\partial c \partial Q$ implies the equality between the relations expressed in Eq. (132) and (67). This equality provides an independent verification of the correctness of the respective expressions as well as a verification of the solution accuracy of computing the second-level adjoint functions $\boldsymbol{\Psi}_1^{(2)} \equiv \left(\Psi_{11}^{(2)}, \Psi_{12}^{(2)}\right)$ and $\boldsymbol{\Psi}_5^{(2)} \equiv \left(\Psi_{51}^{(2)}, \Psi_{52}^{(2)}\right)$. The symmetry of the second-order sensitivity $\partial^2 T(z_r)/\partial c \partial q$ implies the equality between the relations expressed in Eq. (133) and (86). This equality provides an independent verification of the correctness of the respective expressions as well as a verification of the solution accuracy of computing the second-level adjoint functions $\boldsymbol{\Psi}_2^{(2)} \equiv \left(\Psi_{21}^{(2)}, \Psi_{22}^{(2)}\right)$ and $\boldsymbol{\Psi}_5^{(2)} \equiv \left(\Psi_{41}^{(2)}, \Psi_{42}^{(2)}\right)$. Furthermore, the symmetry of the second-order sensitivity $\partial^2 T(z_r)/\partial c \partial T_a$ implies the equality between the relations expressed in Eq. (134) and (104). This equality provides an independent verification of the correctness of the respective expressions as well as a verification of the solution accuracy of computing the second-level adjoint functions $\boldsymbol{\Psi}_3^{(2)} \equiv \left(\Psi_{31}^{(2)}, \Psi_{32}^{(2)}\right)$ and $\boldsymbol{\Psi}_5^{(2)} \equiv \left(\Psi_{41}^{(2)}, \Psi_{42}^{(2)}\right)$. Finally, the symmetry of the second-order sensitivity $\partial^2 T(z_r)/\partial c \partial k_0$ implies the equality between the relations expressed in Eq. (135) and (120). This equality provides an independent verification of the correctness of the respective expressions as well as a verification of the solution accuracy of computing the second-level adjoint functions $\boldsymbol{\Psi}_4^{(2)} \equiv \left(\Psi_{41}^{(2)}, \Psi_{42}^{(2)}\right)$ and $\boldsymbol{\Psi}_5^{(2)} \equiv \left(\Psi_{41}^{(2)}, \Psi_{42}^{(2)}\right)$. The 2nd-order relative sensitivity $\frac{\partial^2 T(z_r)}{\partial c^2} \frac{c^2}{T(z_r)}$ is depicted in Figure 21 as a function of the arbitrary location $z_r \in [-\ell/2, \ell/2]$.

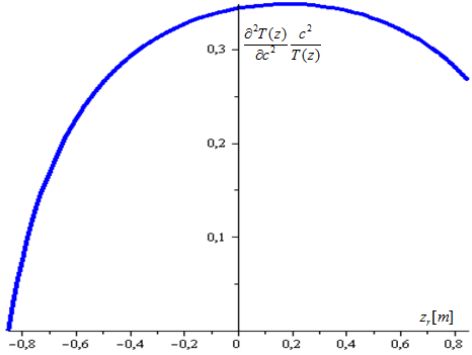

Fig. 21: Relative sensitivity $\dfrac{\partial^2 T(z_r)}{\partial c^2} \dfrac{c^2}{T(z_r)}$

The computation of the 1$^{st}$- and 2$^{nd}$-order sensitivities presented in this Section have underscored that:

(i) *One adjoint computation* was needed to determine the 1$^{st}$-level adjoint function $\Psi(z)$, which sufficed to compute, using just quadratures, all of the first-order sensitivities $S_i \equiv \partial R/\partial \alpha_i$, $i = 1,2,3,4,5$;

(ii) The magnitudes of the 2$^{nd}$-order sensitivities are generally of the same order as the magnitudes of the 1$^{st}$-order sensitivities;

(iii) *Five adjoint computations*, one each for computing the respective 2$^{nd}$-level adjoint functions $\mathbf{\Psi}_i^{(2)} \equiv \left(\Psi_{i1}^{(2)}, \Psi_{i2}^{(2)}\right)$, $i = 1,2,3,4,5$, and subsequently determining all of the 2$^{nd}$-order sensitivities $S_{ij} \equiv \partial^2 R/\partial \alpha_i \partial \alpha_j$; $i,j = 1,2,3,4,5$. These computations also provide independent solution verifications of the 2$^{nd}$-level adjoint functions $\mathbf{\Psi}_i^{(2)} \equiv \left(\Psi_{i1}^{(2)}, \Psi_{i2}^{(2)}\right)$, $i = 1,2,3,4,5$, due to the symmetry of the 2$^{nd}$-order sensitivities.

Notably, similar differential operators appear on the left-sides of the *1$^{st}$- and 2$^{nd}$-LASS*, which therefore would use the same "solvers" of differential equations; *only the sources on the right sides of these suystems differ from each other*. The 1$^{st}$- and 2$^{nd}$-order sensitivities will be used in the following Section to illustrate their essential role for quantifying standard deviations and non-Gaussian features (e.g., asymmetries) of the various response distributions. To quantify assymetries in distribution, at (the very) least the third order ("skewness") response correlations need to be computed, which require the exact computation of (at least) the first- and second-order response sensitivities to model parameters.

## 4. APPLICATION OF THE 2$^{nd}$-ASAM FOR QUANTIFYING NON-GAUSSIAN FEATURES OF THE RESPONSE UNCERTAINTY DISTRIBUTION

In general, the model parameters are experimentally derived quantities and are therefore subject to uncertainties. Specifically, consider that the model comprises $N_\alpha$ uncertain parameters $\alpha_i$, which constitute the components of the (column) vector $\boldsymbol{\alpha}$ of *model parameters*, defined as $\boldsymbol{\alpha} = (\alpha_1, ..., \alpha_{N_\alpha})$. The usual information available in practice comprises the mean values of the model parameters together with uncertainties (standard deviations and, occasionally, correlations) computed about the respective mean values. The components of vector $\boldsymbol{\alpha}^0 \equiv (\alpha_1^0, ..., \alpha_{N_\alpha}^0)$ of mean values of the model parameters are denoted as $\alpha_i^0$ and defined as

$$\alpha_i^0 \equiv \langle \alpha_i \rangle, \quad \langle f \rangle \equiv \int f(\boldsymbol{\alpha}) p(\boldsymbol{\alpha}) d\boldsymbol{\alpha}, . \tag{137}$$

where the angular brackets denotes integration of a generic function $f(\boldsymbol{\alpha})$ over the *unknown* joint probability distribution, $p(\boldsymbol{\alpha})$, of the parameters $\boldsymbol{\alpha}$. The parameter distribution's second-order central moments, $\mu_2^{ij}(\boldsymbol{\alpha})$, are defined as

$$\mu_2^{ij}(\boldsymbol{\alpha}) \equiv \langle (\alpha_i - \alpha_i^0)(\alpha_j - \alpha_j^0) \rangle \equiv \rho_{ij} \sigma_i \sigma_j; \quad i, j = 1, ..., N_\alpha. \tag{138}$$

The central moments $\mu_2^{ii}(\boldsymbol{\alpha}) \equiv \text{var}(\alpha_i)$ are called the *variance* of $\alpha_i$, while the central moments $\mu_2^{ij}(\boldsymbol{\alpha}) \equiv \text{cov}(\alpha_i, \alpha_j); i \neq j$, are called the *covariances* of $\alpha_i$ and $\alpha_j$. The *standard* *deviation* of $\alpha_i$ is defined as $\sigma_i \equiv \sqrt{\mu_2^{ii}(\boldsymbol{\alpha})}$. When the model under consideration is used to compute $N_r$ *responses* (or results), denoted in vector form as $\mathbf{r} = (r_1, ..., r_{N_r})$, each of these responses will be implicit functions of the model's parameters, i.e., $\mathbf{r} = \mathbf{r}(\boldsymbol{\alpha})$. It follows that $\mathbf{r} = \mathbf{r}(\boldsymbol{\alpha})$ will be a vector-valued variate which obeys a (generally intractable) multivariate distribution in $\boldsymbol{\alpha}$. For large-scale systems, the probability distribution $p(\boldsymbol{\alpha})$ is not known in practice and, even if it were known, the induced distribution in $\mathbf{r} = \mathbf{r}(\boldsymbol{\alpha})$ would still be intractable, since $p(\boldsymbol{\alpha})$ could not be propagated exactly through the large-scale models used in for simulating realistic physical systems. The uncertainties in a response $\mathbf{r}(\boldsymbol{\alpha})$ arising from uncertainties in the parameters $\boldsymbol{\alpha}$ can be computed by expanding formally the response $\mathbf{r}(\boldsymbol{\alpha})$ in a Taylor series around $\boldsymbol{\alpha}^0$, constructing appropriate products of such Taylor series, and integrating formally the various products over the unknown parameter distribution function $p(\boldsymbol{\alpha})$, to obtain response correlations. This method for constructing response correlations stemming from parameter correlations is known as the "propagation of errors" or "propagation of moments" method [see, e.g., Ref 7].

For illustrating the effects of second-order response sensitivities for the paradigm nonlinear heat conduction benchmark considered in this work, it suffices to take from Ref. [7] response correlations up to third-order, for the very simple case when: (i) the parameters are uncorrelated and normally distributed; and (ii) only the first- and second-order response sensitivities are available. For these particular conditions, the response correlations derived in [7] reduce to the following expressions for the first three response moments:

(i) The expected value of a response $r_k$, denoted here as $[E(r_k)]^{UG}$, which arises due to uncertainties in uncorrelated normally-distributed model parameters (the superscript *UG* indicates "uncorrelated Gaussian" parameters), is given by the expression

$$[E(r_k)]^{UG} = r_k(\boldsymbol{\alpha}^0) + \frac{1}{2}\sum_{i=1}^{N_\alpha} \frac{\partial^2 r_k}{\partial \alpha_i^2} \sigma_i^2, \tag{139}$$

where $r_k(\boldsymbol{\alpha}^0)$ denotes the computed nominal value of the response;

(ii) The covariance, $cov(r_k, r_\ell)$, between two responses $r_k$ and $r_\ell$ arising from normally-distributed uncorrelated parameters is given by

$$[cov(r_k, r_\ell)]^{UG} = \sum_{i=1}^{N_\alpha} \left(\frac{\partial r_k}{\partial \alpha_i}\frac{\partial r_\ell}{\partial \alpha_i}\right)\sigma_i^2 + \frac{1}{2}\sum_{i=1}^{N_\alpha}\left(\frac{\partial^2 r_k}{\partial \alpha_i^2}\right)\left(\frac{\partial^2 r_\ell}{\partial \alpha_i^2}\right)\sigma_i^4. \tag{140}$$

The variance, $var(r_k)$, of a response $r_k$ is obtained by setting $r_k \equiv r_l$ in the above expression to obtain

$$[var(r_k)]^{UG} = \sum_{i=1}^{N_\alpha} \left(\frac{\partial r_k}{\partial \alpha_i}\right)^2 \sigma_i^2 + \frac{1}{2}\sum_{i=1}^{N_\alpha}\left(\frac{\partial^2 r_k}{\partial \alpha_i^2}\right)^2 \sigma_i^4. \tag{141}$$

As indicated by the expressions in Eqs. (139) - (141), the second-order sensitivities have the following impacts on the response moments:

(a) They cause the "expected value of the response", $[E(r_k)]^{UG}$, to differ from the "computed nominal value of the response", $r_k(\boldsymbol{\alpha}^0)$;

(b) They contribute to the response variances and covariances; however, since the contributions involving the second-order sensitivities are multiplied by the fourth power of the parameters' standard deviations, the total of these contributions is expected to be relatively smaller than the contributions stemming from the first-order response sensitivities.

## 4.1. Computation of Response Standard Deviations

To illustrate the impact of $1^{st}$-order versus $2^{nd}$-order sensitivities in contributing to the uncertainty in the temperature response $T(z_r)$ at an arbitrary location $z_r$, consider that the model parameters $Q$, $q$, $T_a$, $k_0$, $c$ are uncorrelated and normally distributed, all having relative standard deviations of 10%, i.e., $\sigma(Q)/Q^0 = \sigma(q)/q^0 = \sigma(T_a)/T_a^0 = \sigma(k_0)/k_0^0 = \sigma(c)/c^0 = 10\%$, where $\sigma(Q)$, $\sigma(q)$, $\sigma(T_a)$, $\sigma(k_0)$, and $\sigma(c)$ denote the respective absolute standard deviations. It follows from Eq. (141) that the contributions stemming from the $1^{st}$-order parameter derivatives (in the absence of $2^{nd}$- and higher-order derivatives) yields the following "$1^{st}$-order standard deviation," denoted as $\text{Std.Dev}_{FO}[T(z_r)]$, of $T(z_r)$:

$$\text{Std.Dev}_{FO}[T(z_r)] \equiv \left\{ \left[\frac{\partial T(z_r)}{\partial Q}\sigma(Q)\right]^2 + \left[\frac{\partial T(z_r)}{\partial q}\sigma(q)\right]^2 \right. \\ \left. + \left[\frac{\partial T(z_r)}{\partial T_a}\sigma(T_a)\right]^2 + \left[\frac{\partial T(z_r)}{\partial k_0}\sigma(k_0)\right]^2 + \left[\frac{\partial T(z_r)}{\partial c}\sigma(c)\right]^2 \right\}^{1/2}. \tag{142}$$

The magnitudes of each of the terms on the right side of Eq. (142) quantifies the contribution made by each of the model parameters, considerd to by uncorrelated and normally distributed to $1^{st}$-order the stadard deviation $\text{Std.Dev}_{FO}[T(z_r)]$. These magnitudes are plotted, as a function of $z_r$, in Figures 22 – 26.

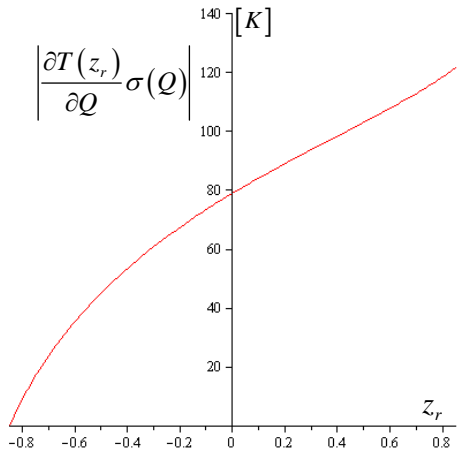

Fig. 22: Contributions to *Std. Dev.* in $T(z_r)$ of 1st-order sensitivities of *Q*.

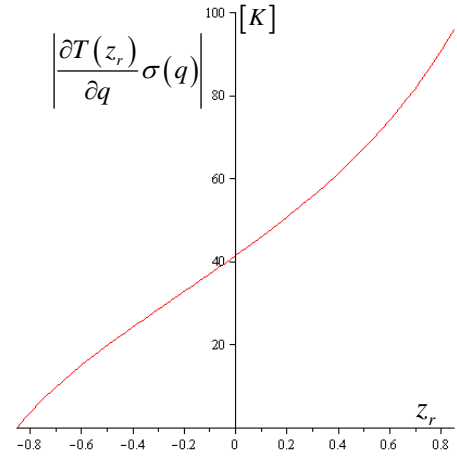

Fig. 23: Contributions to *Std. Dev.* in $T(z_r)$ of 1st-order sensitivities of *q*.

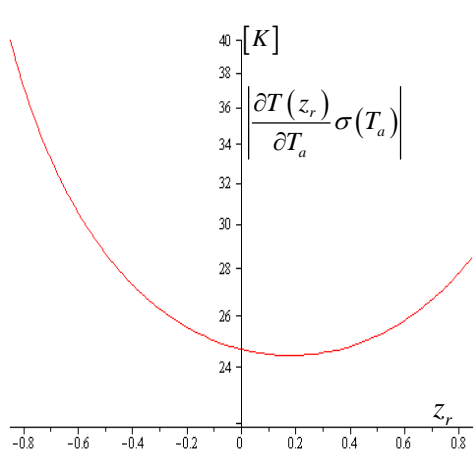

Fig. 24: Contributions to *Std. Dev.* in $T(z_r)$ of 1st-order sensitivities of $T_a$

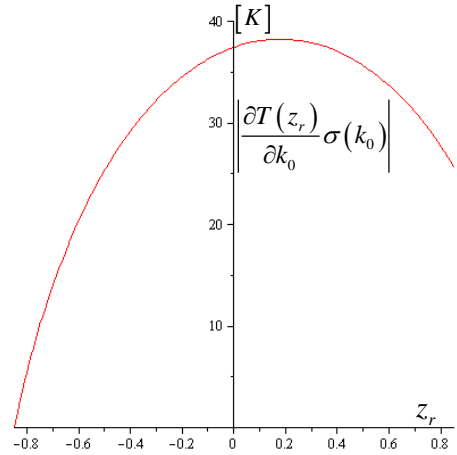

Fig. 25: Contributions to *Std. Dev.* in $T(z_r)$ of 1st-order sensitivities of $k_0$.

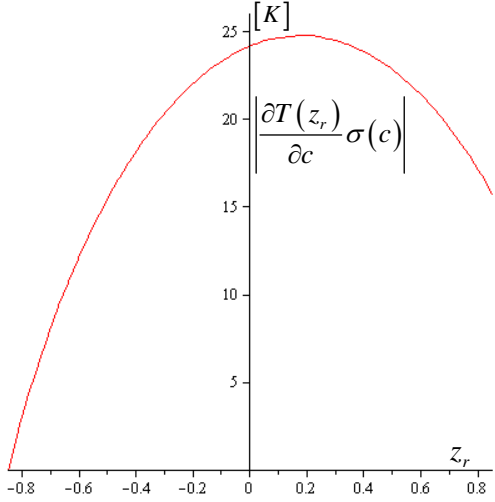
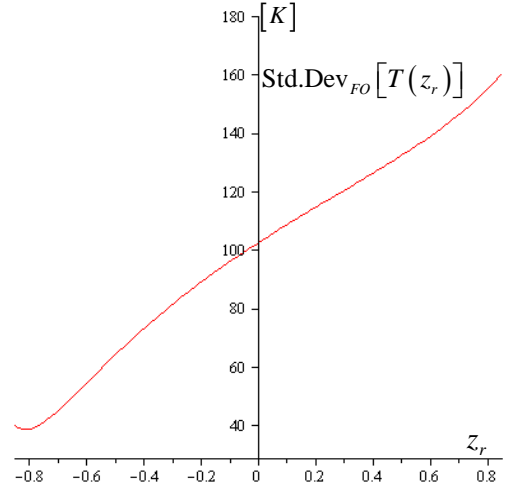

Fig. 26: Contributions to *Std. Dev.* in $T(z_r)$ of 1$^{st}$-order sensitivities of *c*.

Fig. 27: Contributions to *Std. Dev.* in $T(z_r)$ of 1$^{st}$-order sensitivities of in all parameters.

As shown in Figures 22 through 26, the quantities $|\sigma(Q)\partial T(z_r)/\partial Q|$ and $|\sigma(q)\partial T(z_r)/\partial q|$ are monotonically increasing, the quantity $|\sigma(T_a)\partial T(z_r)/\partial T_a|$ displays a minimum, while the quantities $|\sigma(k_0)\partial T(z_r)/\partial k_0|$ and $|\sigma(c)\partial T(z_r)/\partial c|$ display maxima, as functions of $z_r$. Each of these behaviors is governed by the behavior of the respective 1-st order derivatives, of course. Since the derivatives $|\partial T(z_r)/\partial Q|$ and $|\partial T(z_r)/\partial q|$ are the largest, the respective contributions cause the quantity $\text{Std.Dev}_{FO}[T(z_r)]$ to increase monotonically as a function of $z_r$. Noticably, the smallest contribution to $\text{Std.Dev}_{FO}[T(z_r)]$ stems from $|\sigma(c)\partial T(z_r)/\partial c|$, indicating that the parameter *c*, which actually controls the strength of the nonlinearity in the conduction equation, is not very important in the temperature range under consideration (400 – 800 K).

It follows from Eq. (141) that the contributions stemming from the 1$^{st}$-order parameter derivatives (in the absence of 2$^{nd}$- and higher-order derivatives) yields the following "1$^{st}$-order standard deviation," denoted as $\text{Std.Dev}_{FO}[T(z_r)]$, of $T(z_r)$:

The quantities $\left(\sigma_i^2 \partial^2 r_k / \partial \alpha_i^2\right)^2 / 2$ contain the contributions involving the 2$^{nd}$-order derivatives to the total standard deviation, denoted as $\text{Std.Dev}\left[T(z_r)\right]$, of $T(z_r)$. The magnitudes of each of these 2$^{nd}$-order contributions are plotted, as a function of $z_r$, in Figures 28 – 32.

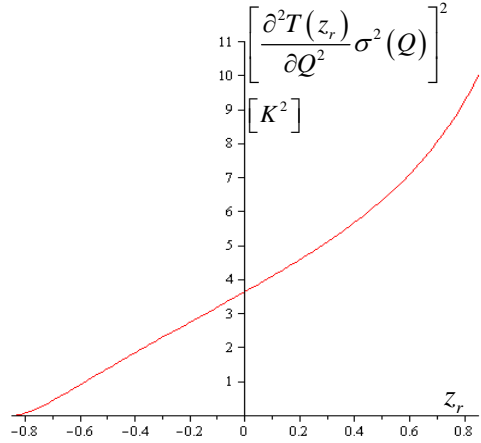

Fig. 28: Contributions from 2$^{nd}$-order sensitivities of $Q$ to *Std. Dev.* in $T(z_r)$.

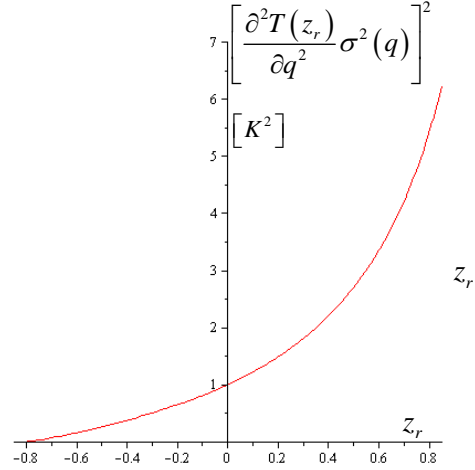

Fig. 29: Contributions from 2$^{nd}$-order sensitivities of $q$ to *Std. Dev.* in $T(z_r)$.

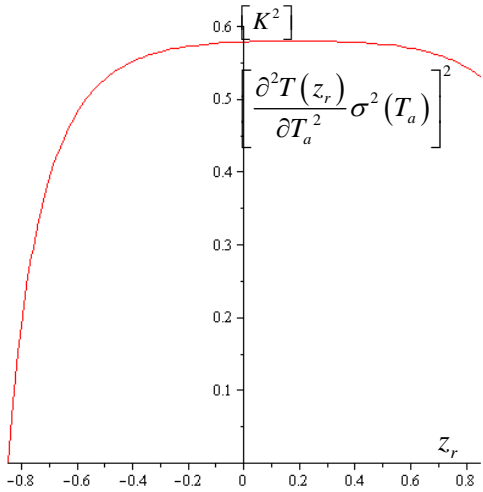

Fig. 30: Contributions from 2$^{nd}$-order sensitivities of $T_a$ to *Std. Dev.* in $T(z_r)$.

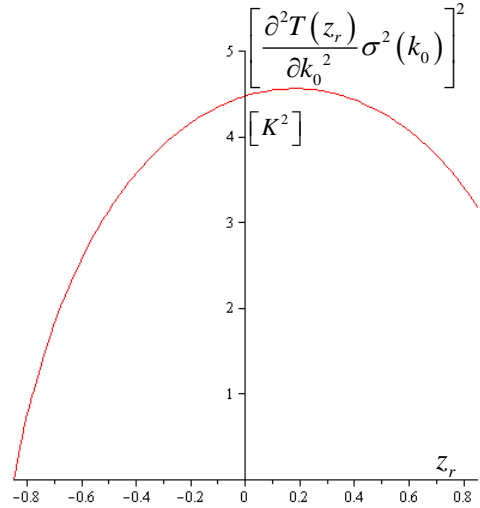

Fig. 31: Contributions from 2$^{nd}$-order sensitivities of $k_0$ to *Std. Dev.* in $T(z_r)$.

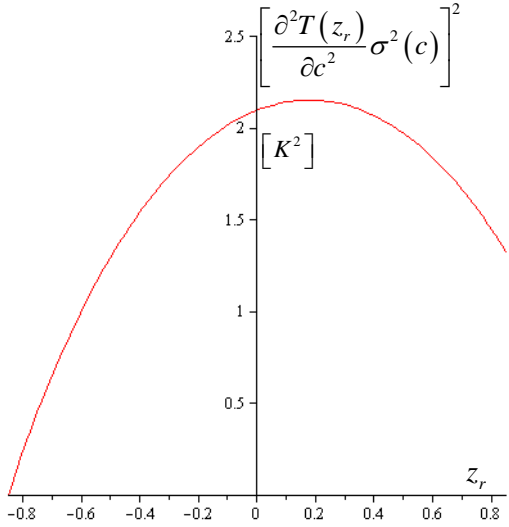 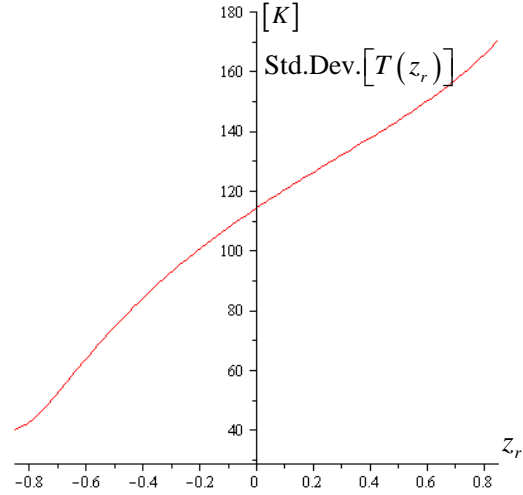

Fig. 32: Contributions of 2nd-order sensitivities of $c$ to *Std. Dev.* in $T(z_r)$.

Fig. 33: Spatial variation of the standard deviation of the temperature distribution $T(z_r)$.

The results plotted in Figures 22-26, to the corresponding results plotted in Figures 28-32 indicate that although some of the 2nd-order sensitivities have magnitudes comparable to the 1st-order ones, the contributions stemming from the 1st-order sensitivities are much smaller than those stemming from the 2nd-order sensitivities. This is because the contributions stemming from the 2nd-order sensitivities are multiplied by the fourth-power of the respective standard deviations, rather than by the second-power, as are the 1st-order sensitivities. The contributions to the total standard deviation in $T(z_r)$ stem predominantly from the heat source $Q$ and the boundary heat flux $q$, which cause this standard deviation to increase monotonically as a function of the location $z_r$, from the inlet to the outlet, reaching its maximum value of 170K at the outlet. Equivalently, the relative standard deviation of $T(z_r)$ increases from 10% the inlet to 24% at the outlet; recall that all of the model parameters were assumed to have a relative standard deviation of 10%. Table 1 presents the actual values at the inlet, outlet, and at $z_r = z_{max} = 18\ cm$, which denotes the location where the temperature distribution $T(z_r)$ reaches its maximum value, $T_{max}(z_{max}) = 874.7K$. Notably, the terms involving the 2nd-order sensitivities do not contribute at all at thebottom of the test section, but their contributions increase monotonically from the bottom to the test section's top.

Table 1: Individual and total contributions of parameters' standard deviations to the 1$^{st}$-order and total standard deviation in $T(z_r)$, at selected locations.

| Location | $z_r = -l/2$ | $z_r = z_{\max}$ | $z_r = l/2$ |
|---|---|---|---|
| Std.Dev$_{Q1}[T(z_r)]$, [K] | 0 | 88 | 122 |
| Std.Dev$_{Q2}[T(z_r)]$, [K] | 0 | 5 | 10 |
| Std.Dev$_{q1}[T(z_r)]$, [K] | 0 | 50 | 96 |
| Std.Dev$_{q2}[T(z_r)]$, [K] | 0 | 1.4 | 6 |
| Std.Dev$_{Ta1}[T(z_r)]$, [K] | 40 | 24 | 29 |
| Std.Dev$_{Ta2}[T(z_r)]$, [K] | 0 | 0.5 | 0.6 |
| Std.Dev$_{k1}[T(z_r)]$, [K] | 0 | 38 | 26 |
| Std.Dev$_{k2}[T(z_r)]$, [K] | 0 | 3 | 5 |
| Std.Dev$_{c1}[T(z_r)]$, [K] | 0 | 25 | 16 |
| Std.Dev$_{c2}[T(z_r)]$, [K] | 0 | 1.3 | 2 |
| Std.Dev$_{FO}[T(z_r)]$, [K] | 40 | 117 | 160 |
| Std.Dev$[T(z_r)]$, [K] | 40 | 125 | 170 |

## 4.2. Computation of Response Skewness

The third-order response correlation, $\mu_3(r_k, r_l, r_m)$, among three responses ($r_k$, $r_\ell$ and $r_m$) has the following expression:

$$\left[\mu_3(r_k, r_l, r_m)\right]^{UG} = \sum_{i=1}^{N_\alpha} \left( \frac{\partial^2 r_k}{\partial \alpha_i^2} \frac{\partial r_l}{\partial \alpha_i} \frac{\partial r_m}{\partial \alpha_i} + \frac{\partial r_k}{\partial \alpha_i} \frac{\partial^2 r_l}{\partial \alpha_i^2} \frac{\partial r_m}{\partial \alpha_i} + \frac{\partial r_k}{\partial \alpha_i} \frac{\partial r_l}{\partial \alpha_i} \frac{\partial^2 r_m}{\partial \alpha_i^2} \right) \sigma_i^4. \qquad (143)$$

In particular, the third-order central moment, $\mu_3(r_k)$, of the response $r_k$ is determined by setting $k = l = m$ in Eq. (143) to obtain

$$\left[\mu_3(r_k)\right]^{UG} = 3 \sum_{i=1}^{N_\alpha} \left( \frac{\partial r_k}{\partial \alpha_i} \right)^2 \frac{\partial^2 r_k}{\partial \alpha_i^2} \sigma_i^4. \qquad (144)$$

The *skewness*, $\gamma_1(r_k)$, of a response $r_k$ can be computed using the customary definition

$$\gamma_1(r_k) \equiv \mu_3(r_k) / [\text{var}(r_k)]^{3/2}. \tag{145}$$

Recall that the *skewness* of a distribution quantifies the departure of the subject distribution from symmetry. Symmetric univariate distributions (e.g., the Gaussian) are characterized by $\gamma_1(r_k) = 0$. A distribution with a long right tail would have a positive skewness while a distribution with a long left tail would have a negative skewness. In other words, if $\gamma_1(r_k) < 0$, then the respective distribution is skewed towards the left of the mean $[E(r_k)]^{UG}$, favoring lower values of $r_k$ relative to $[E(r_k)]^{UG}$. On the other hand, if $\gamma_1(r_k) > 0$, then the respective distribution is skewed towards the right of the mean $[E(r_k)]^{UG}$, favoring higher values of $r_k$ relative to $[E(r_k)]^{UG}$.

As Eq. (143) indicates, neglecting the second-order sensitivities for normally distributed model parameters would nullify the third-order response correlations and hence would nullify the *skewness*, $\gamma_1(r_k)$, of a response $r_k$, cf. Eq. (144). Consequently, any events occurring in a response's long and/or short tails, which are characteristic of rare but important events, would likely be missed. It is of interest to quantify the skewness induced in the temperature distribution $T(z_r)$ by each model parameter considered separately. These "individually-induced" skewnesses in $T(z_r)$ will be denoted as $\gamma_Q[T(z_r)]$, $\gamma_q[T(z_r)]$, $\gamma_{Ta}[T(z_r)]$, $\gamma_k[T(z_r)]$, and $\gamma_c[T(z_r)]$, respectively. In other words, the quantity $\gamma_Q[T(z_r)]$ denotes the skewness that would be induced in the temperature distribution $T(z_r)$ if only the heat source $Q$ were a normally distributed variate with a relative standard deviation of 10%. The quantity $\gamma_q[T(z_r)]$ denotes the skewness that would be induced in the temperature distribution $T(z_r)$ if only the boundary heat flux $q$ were a normally distributed variate with a relative standard deviation of 10%, and so on. The respective 3$^{rd}$-order moments of the temperature distribution $T(z_r)$ are plotted along with the corresponding skewnesses in Figures 34-38, as functions of the location $z_r$, from the bottom to the top of the test section.

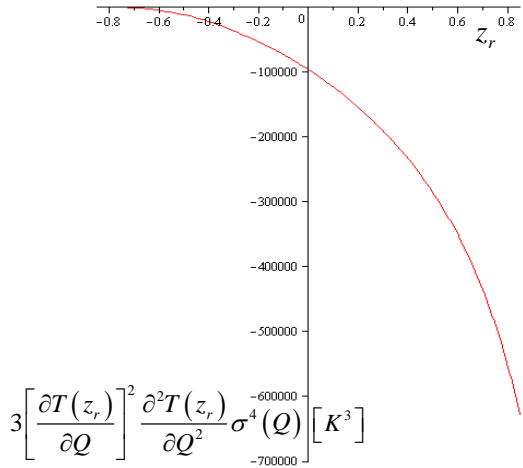

Fig. 34: 3$^{rd}$-order moment of $T(z_r)$ arising solely from the variate $Q$.

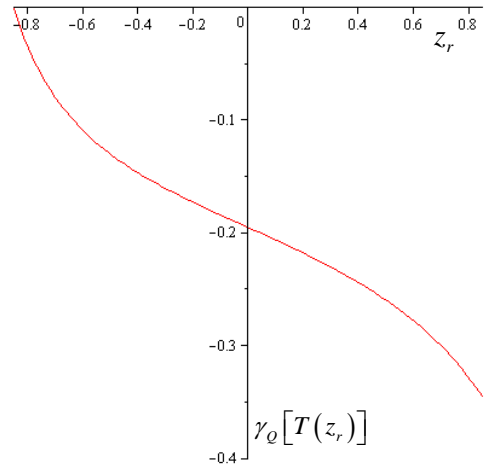

Fig. 35: Skewness of $T(z_r)$ due soley to the variate $Q$.

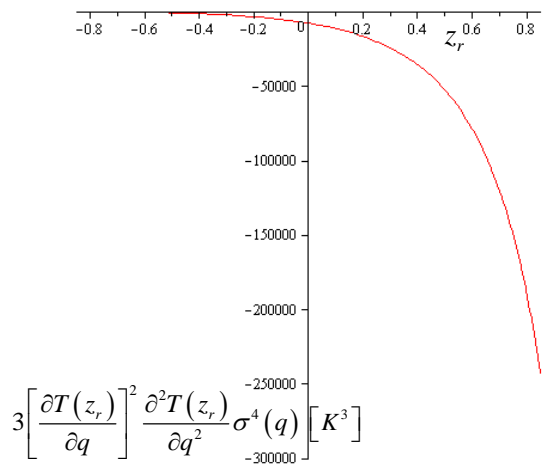

Fig. 36: 3$^{rd}$-order moment of $T(z_r)$ arising solely from the variate $q$.

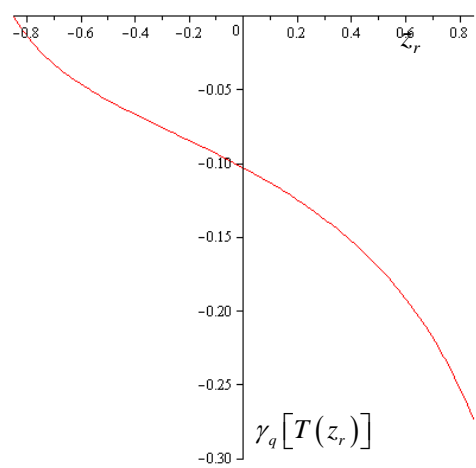

Fig. 37: Contributions to skewness in $T(z_r)$ due solely to the variate $q$.

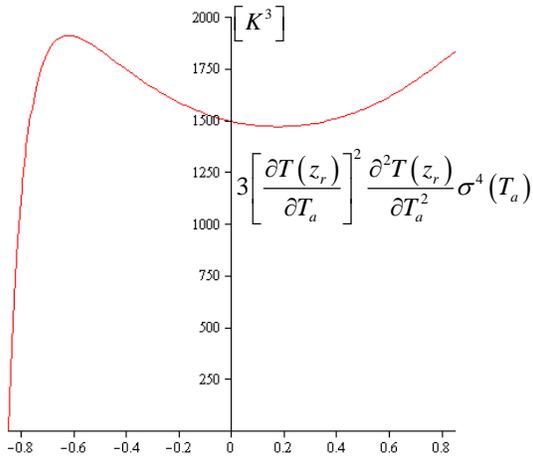

Fig. 38: 3rd-order moment of $T(z_r)$ arising solely from the variate $T_a$.

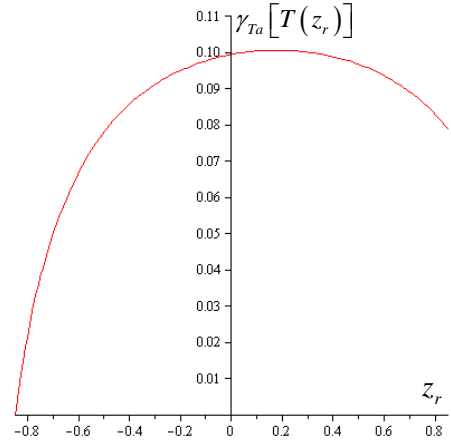

Fig. 39: Contributions to skewness in $T(z_r)$ due solely to the variate $T_a$.

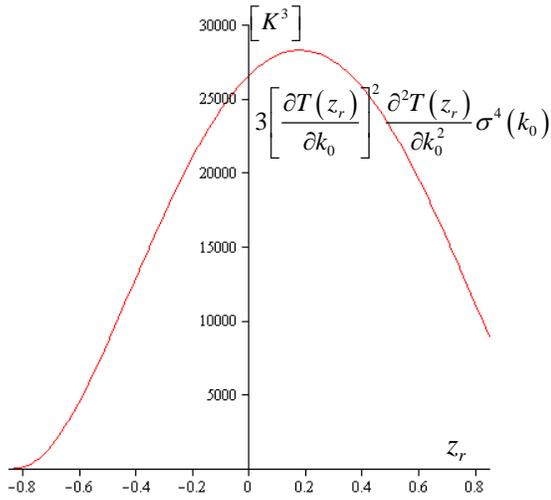

Fig. 40: 3rd-order moment of $T(z_r)$ arising solely from the variate $k_0$.

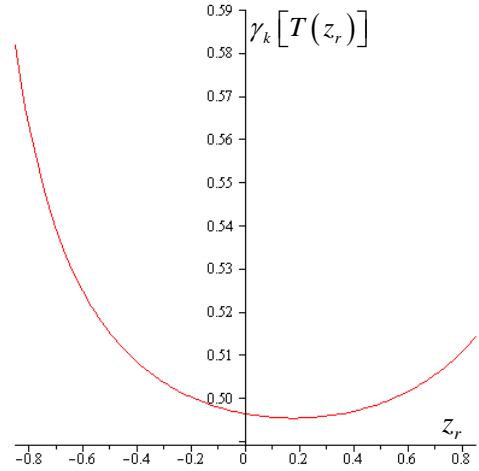

Fig. 41: Contributions to skewness in $T(z_r)$ due solely to the variate $k_0$.

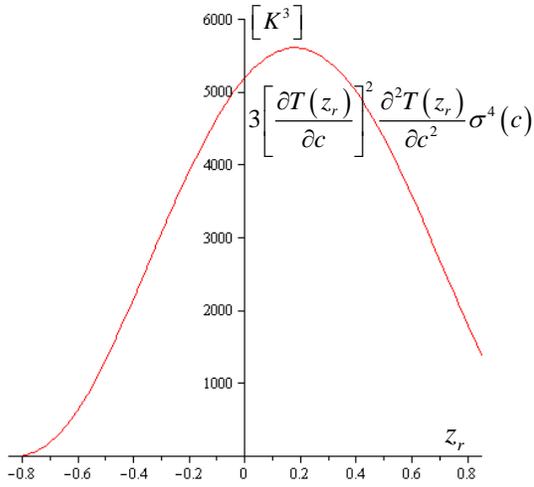
Fig. 42: 3rd-order moment of $T(z_r)$ arising solely from the variate $c$.

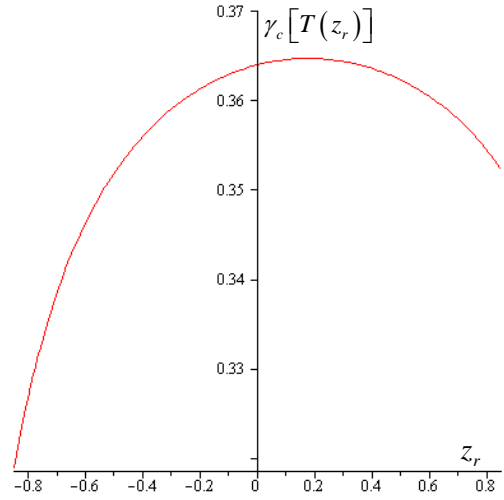
Fig. 43: Contributions to skewness in $T(z_r)$ due solely to the variate $c$.

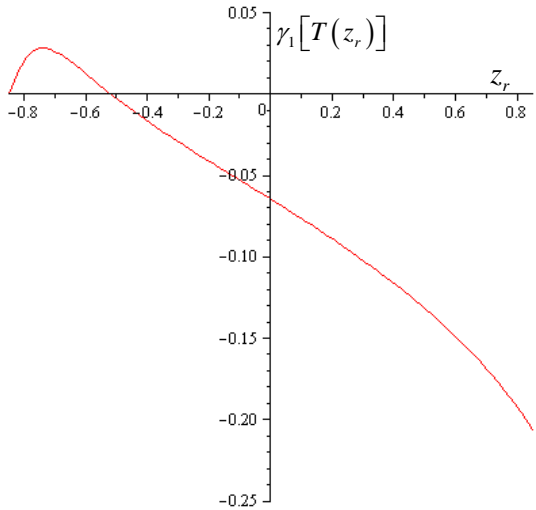
Fig. 44: Spatial variation of the skewness in $T(z_r)$

Figures 35 and 37 indicate that $\gamma_Q[T(z_r)]<0$ and $\gamma_q[T(z_r)]<0$. Hence, if only the heat source $Q$ and/or the boundary heat flux $q$ were normally distributed variates, each having a relative standard deviation of 10%, they would cause the temperature distribution $T(z_r)$ to be skewed significantly towards values *lower* the the mean temperature (i.e., the distribution would be skewed to the "left" of the mean value). On the other hand, Figures 41 and 43 indicate that $\gamma_k[T(z_r)] > 0$ and $\gamma_c[T(z_r)] > 0$. Hence, if only the model parameters $k_0$ and/or $c$ were normally distributed variates, each having a relative standard deviation of 10%, they would cause the temperature distribution $T(z_r)$ to be skewed significantly towards values *higher* the the mean temperature (i.e., the distribution would be skewed to the "right" of the mean value). Finally, Figure 39 indicates that $\gamma_{Ta}[T(z_r)] > 0$, but has rather small positive values. Hence, if the temperature $T_a$ at the bottom of the test section were the sole variate, normally distributed and having a relative standard deviation of 10% (and all other parameters were perfectly known, taking on exactly their nominal values), then the temperature distribution $T(z_r)$ would display a slight asymmetry towards temperatures higher than the mean temperature. The total skewness of the temperature distribution $T(z_r)$, reflecting the cumulative effects of the model parameters (assumed to be normally distributed and having all relative standard deviations of 10%) is depicted in Figure 43. This figure indicates that the negative contributions stemming from the $3^{rd}$-moments of heat source, $Q$, and the boundary flux, $q$, are initially small, so the positive contributions from the $3^{rd}$-order moments of the other model parameters dominate at the bottom of the test section. Therefore, the temperature distribution $T(z_r)$ is skewed slightly toward higher temperatures in the region extending about 20 cm from the bottom of the test section. For the remainder of the test section, however, the negative contributions stemming from the $3^{rd}$-moments of heat source, $Q$, and the boundary flux, $q$, become dominant, so the temperature distribution $T(z_r)$ becomes increasingly skewed towards temperatures lower than (i.e., to the "left" of) the mean temperature towards the top of the test section. Notably, the influemce of the model parameter $c$, which controls the strength of the nonlinearity in this illustrative benchmark problem, would be strong if it were the only uncertain model parameter. However, if all of the other parameters are also uncertain, all having equal relative standard deviations, the uncertainties in the heat source $Q$ and boundary heat flux $q$ estompate the impact uncertainties in $c$, for the range of temperatures (400-900K) considered for this benchmark problem. The numerical values of the impact of the individual parameters, as well as their cumulative impact on the temperature distribution $T(z_r)$ at the bottom and top of the test section, as well as at the location $z_r = z_{max} = 18\ cm$, are presented in Table 2.

Table 2: Individual and total contributions of the parameters' 2$^{nd}$-order sensitivities to the skewness, $\gamma_1[T(z_r)]$, in $T(z_r)$, at selected locations.

| Location | $z_r = -l/2$ | $z_r = z_{max}$ | $z_r = l/2$ |
|---|---|---|---|
| $\gamma_Q[T(z_r)]$ | 0 | -0.22 | -0.35 |
| $\gamma_q[T(z_r)]$ | 0 | -0.12 | -0.27 |
| $\gamma_{Ta}[T(z_r)]$ | 0 | 0.1 | 0.08 |
| $\gamma_k[T(z_r)]$ | 0 | 0.49 | 0.51 |
| $\gamma_c[T(z_r)]$ | 0 | 0.36 | 0.35 |
| $\gamma_1[T(z_r)]$ | 0 | -0.08 | -0.2 |

## 5. SUMMARY AND CONCLUSIONS

This work has presented an illustrative application of the *second-order adjoint sensitivity analysis methodology* (*2$^{nd}$-ASAM*) *for nonlinear systems* to a paradigm nonlinear heat conduction benchmark problem that simulates the steadystate temperature distribution in a conceptual design of a test section comprising heated rods. This benchmark problem is sufficiently simple to admit an exact solution, thereby making transparent the mathematical derivations presented in [1]. The general theory [1] underlying the *2$^{nd}$-ASAM* indicates that, for a physical system comprising $N_\alpha$ parameters, the computation of all of the first- and second-order response sensitivities requires at most $N_\alpha$ "large-scale" computations involving correspondingly constructed adjoint systems, which we call the *second-level adjoint sensitivity systems* (*2$^{nd}$-LASS*). For the illustrative nonlinear heat conduction benchmark considered in this work, six "large-scale" adjoint computations sufficed for the complete and exact computations of all 5 first- and 15 distinct second-order derivatives. It has also been shown that the construction and solution of the *2$^{nd}$-LASS* requires very little additional effort beyond the construction of the adjoint sensitivity system needed for computing the first-order sensitivities. Very significantly, only the sources on the right-sides of the "solver" for the nonlinear heat conduction equation needed to be modified, while the left-side of the differential equation remained unchanged.

For the nonlinear heat conduction benchmark, several $2^{nd}$-order relative sensitivities had magnitudes comparable to those of the $1^{st}$-order relative sensitivities. We have subsequently showed that the second-order sensitivities play the following important roles:

(a) They cause the "expected value of the response" to differ from the "computed nominal value of the response," but for the nonlinear heat conduction benchmark these differences were insignificant.

(b) They contribute to increasing the response variances and modifying the response covariances, but for the nonlinear heat conduction benchmark their contribution is much smaller than (a few % of) that stemming from the first-order response sensitivities.

(c) They comprise the leading contributions to causing asymmetries in the response distribution. For the benchmark test section considered in this work, the heat source, the boundary heat flux, and the temperature at the bottom boundary of the test section would cause the temperature distribution in the test section to be skewed significantly towards values *lower* the the mean temperature. On the other hand, the model parameters entering the nonlinear, temperature-dependent, expression of the LBE conductivity would cause the temperature distribution in the test section to be skewed significantly towards values *higher* the the mean temperature. These opposite effects partially cancel each other. Consequently, the cumulative effects of model parameter uncertainties on the skewness of the temperature distribution in the test section is such that the temperature distribution in the LBE is skewed slightly toward higher temperatures in the cooler part of the test section. However, the temperature distribution becomes increasingly skewed towards temperatures lower than the mean temperature in the hotter part of the test section. Notably, the influence of the model parameter that controls the strength of the nonlinearity in the heat conduction coefficient for this LBE test section benchmark would be strong if it were the only uncertain model parameter. However, if all of the other model parameters are also uncertain, all having equal relative standard deviations, the uncertainties in the heat source and boundary heat flux estompate the impact of uncertainties in the nonlinear heat conduction coefficient for the range of temperatures (400-900K) considered for this LBE test section benchmark.

Ongoing work aims at generalizing the *2$^{nd}$-ASAM* to enable the exact computation of higher-order response sensitivities in an efficient manner, which is expected to affect significantly the fields of optimization and predictive modeling, including uncertainty quantification, data assimilation, model calibration and extrapolation.

**ACKNOWLEDGMENTS**

This work was partially supported by contract 15540-FC51 from Gen4Energy, Inc., with the University of South Carolina.